\documentclass[reqno,11pt]{amsart}
\usepackage[top=2.6cm,left=2.7cm,right=2.7cm,bottom=2.6cm]{geometry}
\usepackage{amsmath,amssymb,amsfonts,amsthm,mathrsfs}
\usepackage[dvipsnames]{xcolor}
\usepackage{tikz}
\usepackage{hyperref}
\hypersetup{
  colorlinks=true,
  citecolor=NavyBlue,
  linkcolor=Brown,
  urlcolor=Orange
}
\newtheorem{conjecture}{Conjecture}[section]
\newtheorem{question}[conjecture]{Question}
\newtheorem{proposition}[conjecture]{Proposition}
\newtheorem{lemma}[conjecture]{Lemma}
\newtheorem{theorem}[conjecture]{Theorem}
\newtheorem{corollary}[conjecture]{Corollary}
\theoremstyle{definition}
\newtheorem{definition}[conjecture]{Definition}
\newtheorem{remark}[conjecture]{Remark}
\newtheorem{example}[conjecture]{Example}
\DeclareMathOperator{\Rep}{Rep}
\DeclareMathOperator{\Lie}{Lie}
\DeclareMathOperator{\Aut}{Aut}
\DeclareMathOperator{\End}{End}
\DeclareMathOperator{\gr}{gr}
\DeclareMathOperator{\MT}{MT}
\newcommand{\Q}{\mathbb Q}
\newcommand{\Z}{\mathbb Z}
\newcommand{\R}{\mathbb R}
\newcommand{\C}{\mathbb C}
\newcommand{\Af}{\mathbb A_{\Q,f}}
\newcommand{\Zhat}{\widehat{\Z}}
\newcommand{\Gm}{\mathbb G_m}
\newcommand{\MM}{\mathcal M}
\newcommand{\Mf}{\mathsf M}
\newcommand{\Fil}{\operatorname{Fil}}
\title[A motivic Vojta inequality for non-CM elliptic curves]
{A motivic Vojta inequality for non-CM elliptic curves with unbounded height}

\author{Daiki Kawabe}
\address{
  Department of Mathematics, School of Education,
  Waseda University,
  Shinjuku, Tokyo, 169-8050, Japan
}
\email{kawabe@aoni.waseda.jp, daiki.kawabe.math@gmail.com}

\date{November 4, 2026}

\begin{document}

\begin{abstract}
We study a bounded-degree form of a motivic Vojta inequality
for pure Andr\'e motives, arising from a question of Kato.
We verify this inequality for infinitely many non-CM elliptic curves
with unbounded stable Faltings height.
The point is that ramification at the finite poles of the
$j$-invariant gives the required discriminant bound.
\end{abstract}
\maketitle

\section{Introduction}
In this paper, we study a bounded-degree form of Kato's motivic Vojta question, using the untruncated counting term and working with pure Andr\'e motives of fixed tensor type and principal level.
Our main result verifies the inequality for infinitely many non-CM elliptic curves with unbounded stable Faltings height and uniformly bounded full-level degree.
We first recall the classical conjecture, then formulate its motivic analogue and present our results.

\subsection{The classical Vojta conjecture}
Fix a number field $K\subset\overline K\subset\C$, a finite set $S$ of
places of $K$ containing the archimedean places, and an integer $r\geq1$.
For a finite extension $F/K$ in $\overline K$, write $S_F$ for the places
above $S$ and $D_F$ for the absolute discriminant.
We use the $K$-relative normalization
\[
 d_K(F):=\frac{\log|D_F|}{[F:K]}-\log|D_K|\geq0,
\]
as in \cite[Definition~24.1]{Vojta}.
Let $X$ be a smooth projective variety over $K$, $D$ a reduced normal
crossings divisor, and $A$ an ample divisor.
Write $K_X$ for a canonical divisor and $h_L$ for
$[K:\Q]$ times an absolute logarithmic Weil height associated
with a divisor $L$ on $X$.
For $x\in X(\overline K)\setminus\operatorname{Supp}(D)$, put $F=K(x)$
and choose compatible local Weil functions $\lambda_{D,v}$.
The counting term is
\[
 N_S(D,x)=\frac1{[F:K]}
 \sum_{v\notin S_F}\lambda_{D,v}(x).
\]

We recall a bounded-degree form of Vojta's conjecture
\cite[Conjecture~25.1 (b)]{Vojta}, allowing the exceptional set
to depend on $r$ as in \cite[Theorem~1.2]{Levin}.
\begin{conjecture}[Vojta, bounded-degree form]
\label{conj:intro-classical}
For every $\varepsilon>0$, there exist a proper closed subset
$Z\subset X$ and a constant $C_\varepsilon\in\R$ such that
\begin{equation}\label{eq:classical}
 h_{K_X+D}(x)\leq d_K(K(x))+N_S(D,x)
                  +\varepsilon h_A(x)+C_\varepsilon
\end{equation}
for every $x\in X(\overline K)\setminus Z$ with
$[K(x):K]\leq r$ and $x\notin\operatorname{Supp}(D)$.
\end{conjecture}

This bounds the log-canonical height by the discriminant,
the counting term, and an arbitrarily small multiple of an ample height.

\subsection{The motivic Vojta inequality and the main results}
We adapt Kato's height and degeneration invariants
\cite[\S\S0.1--0.3]{KatoI} to pure Andr\'e motives within his
group-relative framework \cite[\S\S1.3, 3.2, and 3.6]{KatoII}.
Fixing a pure tensor type $\Phi$, a lattice and principal level,
and a positive height datum $\Lambda$, we associate to an Andr\'e-motivic
object $\mathcal M/F$ a canonical-type height
$h_{\spadesuit,\Phi,K}(\mathcal M)$, a degeneration term
$N^\heartsuit_{K,S}(\mathcal M)$, and a positive height
$h_{\Lambda,K}(\mathcal M)$.
These invariants are defined from the motive and its realizations,
without requiring a moduli space of motives.

Our main focus is the following bounded-degree, untruncated variant of
Kato's motivic Vojta question
\cite[\S2.5.2, Question~4 (1)]{KatoI}.
\begin{question}[A pure-motive variant of Kato's question]
\label{conj:intro-motivic}
Fix $K,S,r,\Phi$, the lattice and principal level, and a
positive height datum $\Lambda$. For which such data does the following
bound hold? For every $\varepsilon>0$, there is
a real constant $C_\varepsilon$, depending only on these fixed data
and $\varepsilon$, such that
\begin{equation}\label{eq:intro-main}
 h_{\spadesuit,\Phi,K}(\mathcal M)
 \leq d_K(F)+N^\heartsuit_{K,S}(\mathcal M)
        +\varepsilon h_{\Lambda,K}(\mathcal M)+C_\varepsilon
\end{equation}
for every motivically generic object $\mathcal M$ of the fixed type
over a finite extension $F/K$ with $[F:K]\leq r$.
\end{question}

The correspondence with \eqref{eq:classical} is
\[
\begin{array}{c|ccccc}
\text{Vojta's conjecture}
& x/K(x)
& h_{K_X+D}(x)
& d_K(K(x))
& N_S(D,x)
& h_A(x)
\\ \hline
\text{Motivic analogue}
& \mathcal M/F
& h_{\spadesuit,\Phi,K}(\mathcal M)
& d_K(F)
& N^\heartsuit_{K,S}(\mathcal M)
& h_{\Lambda,K}(\mathcal M)
\end{array}
\]
We state the question precisely in Question~\ref{conj:main} and compare it with Kato's formulation.

Our first results verify \eqref{eq:intro-main} in cases where the
height terms vanish or are bounded: motives of a single Hodge type,
CM abelian varieties of a fixed tensor type with $G^0$ a torus
and a fixed integral level, and non-CM
elliptic curves of bounded stable Faltings height
(Propositions~\ref{prop:examples-diagonal},
\ref{prop:examples-cm}, and~\ref{prop:examples-bounded-elliptic}).

Our main result verifies \eqref{eq:intro-main} for infinitely many non-CM elliptic curves with unbounded stable Faltings height
(Theorem~\ref{ex:examples-unbounded}).
To explain the main point, let $E/F$ be a non-CM elliptic curve and put
$b_K(E):=[K:\Q]h_{\mathrm{Fal}}^{\mathrm{st}}(E).$
If $E[n]\subset E(F)$, choose a natural full level-$n$ object $\mathcal M_E$
attached to $h^1(E)$. Then
the elliptic height formulas give
\[
h_{\spadesuit,\Phi_{\mathrm{ell}},K}(\mathcal M_E)
-\varepsilon h_{\Lambda,K}(\mathcal M_E)
=(2-\varepsilon\beta_\Lambda)b_K(E),
\qquad \beta_\Lambda>0;
\]
see \eqref{eq:examples-abelian-height}--\eqref{eq:examples-elliptic-height}.
Thus, for $0<\varepsilon<2/\beta_\Lambda$, the left-hand side is
unbounded as $b_K(E)\to+\infty$, so the bounded-height argument no
longer applies.

Instead, we instead use the discriminant $d_K(F)$ to bound $b_K(E)$.
Fix $n\geq3$. For $d\geq2$, put $a_d=1728d+1$ and consider
the elliptic curve
\[
E_{a_d}:\quad y^2=x^3-3a_dx-2a_d^2,
\qquad
j(E_{a_d})=-\frac1d.
\]
If $d$ is cube-free and $F/K$ is a finite extension satisfying
$E_{a_d}[n]\subset E_{a_d}(F)$, then
\[
2b_K(E_{a_d})\leq d_K(F)+C_K,
\]
where $C_K$ is independent of $d$, $F$, and $n$.
The point is that ramification at the finite poles of $j(E_{a_d})$
gives the required lower bound for $d_K(F)$; see
Proposition~\ref{prop:examples-projective-level} and
Proposition~\ref{prop:examples-rational-j}.
Since $N^\heartsuit_{K,S}\geq0$ and $h_{\Lambda,K}$ is bounded below, the bound proves \eqref{eq:intro-main}.

Moreover, as $d$ runs through the primes, these curves give infinitely
many geometric isomorphism classes with unbounded stable Faltings height,
while the full-level fields $F_d=K(E_{a_d}[n])$ have uniformly bounded
degree over $K$.
Theorem~\ref{ex:examples-unbounded} gives a more general statement
for rational $j$-invariants of uniformly bounded absolute value
satisfying \eqref{eq:examples-rational-denominator};
Example~\ref{ex:examples-weierstrass} treats $E_{a_p}$ for primes $p$.

The paper is organized as follows.
Section~2 introduces the motivic framework and the arithmetic invariants
corresponding to the heights and counting terms in Vojta's conjecture.
Section~3 formulates the motivic Vojta inequality, verifies it in
several cases, and proves our main result.

The following diagram summarizes the main dependencies leading to
the unbounded-height result (\S3.9) and its extension to symmetric powers (\S3.10).

\begin{center}
\small
\S2 $\longrightarrow$ \S3.4 $\longrightarrow$ \S3.6
$\longrightarrow$ \S3.8 $\longrightarrow$ \S3.9
$\longrightarrow$ \S3.10
\end{center}

\section{Motivic framework and arithmetic invariants}
\label{sec:motivic-framework}

\subsection{Motives with fixed tensor and integral data}

To define these invariants, we adapt Kato's $G$-motivic framework
to Andr\'e motives.
For each number field $F\subset\C$, let $\MM_F$ be the neutral
semisimple Tannakian category of Andr\'e motives with rational coefficients,
allowing finite direct sums of different pure weights \cite[Theorem~0.4 and \S4.6]{Andre}.
It has the required cohomological realizations without
assuming the standard conjectures.

We fix the data in two stages.
First, we fix a rational tensor type, including the group, the Tate character, the polarization, and the Hodge structures.
We then add a fixed integral structure and level at the finite places.

\subsubsection{A fixed rational tensor type}

Fix a tuple
\[
 \Phi=(G,w_G,\eta,w_0,V_0,q_0,\Upsilon)
\]
consisting of the following data.
\begin{enumerate}
\item \emph{(The group)}
The group $G$ is a linear algebraic group over $\Q$ whose identity
component $G^0$ is reductive. We write $\Rep_\Q(G)$ for its category of
finite-dimensional rational representations.

\item \emph{(The weight cocharacter)}
The map $w_G:\Gm\to Z(G)$ is a central cocharacter.
A representation has weight $w\in\Z$ if $w_G(t)$ acts on it as $t^w$.

\item \emph{(The Tate character)}
The map $\eta:G\to\Gm$ is a character satisfying
$\eta(w_G(t))=t^2$.
Let $\Q_\eta$ denote the representation on $\Q(-1)_B$
on which $G$ acts via $\eta$. Thus $\Q_\eta$ has weight $2$.

\item \emph{(The representation and its weight)}
Let $w_0\in\Z$, and let $V_0\in\Rep_\Q(G)$ be a nonzero representation
of weight $w_0$.
Assume that $V_0\oplus\Q_\eta$ is faithful.
Thus $V_0$ and $\Q_\eta$ generate $\Rep_\Q(G)$
as a Tannakian category.

\item \emph{(The pairing)}
We fix a nondegenerate $G$-equivariant pairing
\[
 q_0:V_0\otimes V_0\rightarrow\Q_\eta^{\otimes w_0},
\]
symmetric for even $w_0$ and alternating for odd $w_0$; negative tensor
powers denote powers of the dual. In particular,
$q_0(gx,gy)=\eta(g)^{w_0}q_0(x,y)$.
The normalization of $\eta$ is compatible with this pairing:
both $V_0\otimes V_0$ and $\Q_\eta^{\otimes w_0}$ have weight $2w_0$.

\item \emph{(The Hodge structures)}
The datum $\Upsilon$ is a $G(\R)$-conjugacy class of homomorphisms
\[
 h:\mathbb S:=\operatorname{Res}_{\C/\R}\Gm\rightarrow G_\R
\]
whose restriction to $\mathbb G_{m,\R}\subset\mathbb S$ is $w_G$.
For $V\in\Rep_\Q(G)$, we use Kato's convention: $V^{p,q}\subset V_\C$
is the subspace on which $h(z)$ acts as $z^p\bar z^q$, and
\[
 \Fil^aV_\C=\bigoplus_{p\geq a}V^{p,q}.
\]
For every $h\in\Upsilon$, we have $\eta(h(z))=z\bar z$,
so $\Q_\eta$ has Hodge type $(1,1)$.
Assume that $q_0$ polarizes the Hodge structure on $V_0$
for some $h\in\Upsilon$.
\end{enumerate}

This is the pure case of the data in
\cite[\S\S1.1.4--1.1.8 and 1.2.9--1.2.10]{KatoII}.

Let $D_\Phi=D(G,\Upsilon)$ be the period domain, which in this pure
setting is $\Upsilon$ with its natural complex manifold structure
\cite[\S\S1.1.8 and 1.1.14--1.1.17]{KatoII}.
Put $\mathfrak g=\Lie(G)$, and for $h\in D_\Phi$ let
$\Fil_h^\bullet\mathfrak g_\C$ be the Hodge filtration induced by
$\operatorname{Ad}\circ h$.
Write $T_{D_\Phi,h}$ for the holomorphic tangent space at $h$
and $T_{D_\Phi,\mathrm{hor},h}$ for its horizontal subspace,
defined by Griffiths transversality
\cite[\S1.2.23]{KatoII}. Then
\begin{equation}\label{eq:tangent}
 T_{D_\Phi,h}\simeq\mathfrak g_\C/\Fil_h^0\mathfrak g_\C,
 \qquad
 T_{D_\Phi,\mathrm{hor},h}\simeq
 \Fil_h^{-1}\mathfrak g_\C/\Fil_h^0\mathfrak g_\C.
\end{equation}

\subsubsection{A fixed integral structure and level}

We fix an integral structure and a principal congruence level.

Put $W=V_0\oplus\Q_\eta$.
Fix a lattice $L_0\subset V_0$, a nonzero vector
$e_\eta\in\Q_\eta$, and an integer $n\geq3$, and put
\[
L=L_0\oplus\Z e_\eta\subset W.
\]
Write $\Af$ for the finite adele ring of $\Q$ and
$\Zhat=\prod_\ell\Z_\ell$, where $\ell$ runs over the primes.
Put $\widehat L=\Zhat\otimes_\Z L$, and set
\begin{equation}\label{eq:level}
\mathcal K=
\{g\in G(\Af):
g\widehat L=\widehat L,\ 
g=1\text{ on }\widehat L/n\widehat L\}.
\end{equation}
The condition $n\geq3$ ensures that $\mathcal K$ satisfies
Kato's neat condition \cite[Remark~1.2.12]{KatoII}.
Throughout, $L_0$, $e_\eta$, and $n$ remain fixed.

\begin{definition}[A motive with fixed rational tensor type and integral level]

An object of type $(\Phi,\mathcal K)$ over $F$ is a pair
$\mathcal M=(\Mf,\lambda),$
where
\[
\Mf:\Rep_\Q(G)\rightarrow\MM_F
\]
is an exact tensor functor sending representations of weight $w$
to pure motives of weight $w$,
and $\lambda=\widetilde\lambda\mathcal K$ is the orbit of an adelic
tensor isomorphism
\[
\widetilde\lambda:
(V\mapsto\Af\otimes_\Q V)
\xrightarrow{\sim}
(V\mapsto\Mf(V)_{\mathrm{et}}).
\]
Here $(-)_{\mathrm{et}}$ denotes the finite adelic realization.
We require this pair to satisfy the following conditions.
\begin{enumerate}
\item \emph{(Galois invariance)}
The orbit $\lambda$ is Galois invariant, meaning that
\[
 \widetilde\lambda^{-1}\circ\sigma\circ\widetilde\lambda
 \in\mathcal K
 \qquad
 \text{for every }\sigma\in\operatorname{Gal}(\overline F/F).
\]

\item \emph{(Betti, Hodge, and Tate compatibility)}
There exists a rational tensor isomorphism
\[
\theta:
(V\mapsto\Mf(V)_B)
\xrightarrow{\sim}
(V\mapsto V)
\]
such that, for some $h\in\Upsilon$, every $\theta_V$
identifies the Hodge structure on $\Mf(V)_B$
with that on $V$ defined by $h$.
Its component $\theta_{\Q_\eta}$ is induced by an isomorphism
$\Mf(\Q_\eta)\simeq\Q(-1)$ in $\MM_F$.
Via this isomorphism, $\Mf(q_0)$ must be a polarization of
$\Mf(V_0)$ or $(-1)^{w_0}$ times a polarization, as in
\cite[\S1.3.3]{KatoII}. The choice of $\theta$ is not part of the data.
\end{enumerate}

Two objects $(\Mf_1,\lambda_1)$ and $(\Mf_2,\lambda_2)$ are
isomorphic if there is a tensor isomorphism
$\alpha:\Mf_1\xrightarrow{\sim}\Mf_2$ over $F$ such that
$
\alpha_{\mathrm{et}}\circ\widetilde\lambda_1
=\widetilde\lambda_2\circ k$ for some $k\in\mathcal K$,
where $\widetilde\lambda_i$ represents $\lambda_i$.
We write
\[
\mathscr M_{\Phi,\mathcal K}(F)
\]
for the set of isomorphism classes of objects of type
$(\Phi,\mathcal K)$ over $F$.

\end{definition}

For such an object, put $M_0=\Mf(V_0)$ and write
$A_{\mathcal M}=\Mf(\mathfrak g)$ for its adjoint motive.

Let $\mathcal C_{\mathcal M}$ be the full Tannakian subcategory generated
by the image of $\Mf$.
Because $W$ is faithful, it is generated by $M_0$ and $\Q(-1)$.
Let $\omega_B$ denote the Betti realization functor. The Betti
Tannakian group of $\mathcal C_{\mathcal M}$ is
\[
 G_{\mathcal M}=\Aut^\otimes
       (\omega_B|_{\mathcal C_{\mathcal M}}).
\]
Choosing $\theta$ as in the preceding definition, the tensor functor
induces a closed immersion $G_{\mathcal M}\hookrightarrow G$.

\begin{definition}[Motivic genericity]

We call $\mathcal M$ \emph{motivic generic} if
\[
 G_{\mathcal M}^0=G^0,
 \qquad\text{equivalently}\qquad
 \dim G_{\mathcal M}=\dim G.
\]

\end{definition}

This is independent of $\theta$.
Equivalently, $\Mf$ does not factor through $\Rep(G')$ for any
algebraic subgroup $G'\subset G$ with $\dim G'<\dim G$
\cite[\S1.3.13]{KatoII}.
If $\MT(\Mf(W)_B)=G^0$, then $\mathcal M$ is motivically generic,
since $\MT(\Mf(W)_B)\subseteq G_{\mathcal M}^0\subseteq G^0$.

\subsection{Canonical-type and positive heights}
We first recall Kato's basic heights of motives, then use them to define the canonical-type and positive heights.
Venkatesh independently studied heights of motives \cite[\S1]{KatoI}
and defined an automorphic height, conjecturing its relation to the
Kato--Koshikawa height \cite{Venkatesh2024}.

\subsubsection{Kato's basic heights}
Let $P/F$ be a pure Andr\'e motive of weight $w$ with a Galois-stable
$\Zhat$-lattice $T(P)\subset P_{\mathrm{et}}$. The pair $(P,T(P))$ is the analogue in the Andr\'e
category of Kato's motive with $\Z$-coefficients
\cite[\S1.1.6]{KatoI}.
Its de Rham realization $P_{\mathrm{dR}}$ is a
finite-dimensional $F$-vector space with a decreasing Hodge filtration
$\Fil^\bullet$; for $P=h^r(X)$ with $X/F$ smooth and projective, it is
$H^r_{\mathrm{dR}}(X/F)$.
Write
$\gr^j_{\Fil}P_{\mathrm{dR}}=\Fil^jP_{\mathrm{dR}}/\Fil^{j+1}P_{\mathrm{dR}}$
and put
\[
 L_j(P)=\det_F\gr_{\Fil}^jP_{\mathrm{dR}}
       \otimes_F(\det_F\gr_{\Fil}^{w-j}P_{\mathrm{dR}})^{-1}.
\]

We give this line Kato's local metrics as follows.
At a finite place $v\mid p$, let $\C_v=\widehat{\overline{F_v}}$,
with valuation ring $\mathcal O_{\C_v}$.
The integral $p$-adic construction applied to the $p$-component
$T(P)_p$ gives lattices on the Hodge graded pieces, hence an
$\mathcal O_{\C_v}$-lattice
\[
\mathcal L_{j,v}\subset\C_v\otimes_F L_j(P)
\]
\cite[\S\S1.3.4--1.3.5 and 1.4.3]{KatoI}.
If $e$ generates this lattice, define
\[
\|a e\|_{j,v}=|a|_v\qquad(a\in\C_v),
\]
where $|\cdot|_v$ extends $|a|_v=(\#k_v)^{-\operatorname{ord}_v(a)}$
on $F_v^\times$, with $k_v$ the residue field and $\operatorname{ord}_v$ the normalized valuation.

At an archimedean place $v$ induced by $\sigma:F\hookrightarrow\C$, put
$D_{j,\sigma}=\C\otimes_{F,\sigma}\det_F\gr^j_{\Fil}P_{\mathrm{dR}}$.
Via Betti--de Rham comparison, complex conjugation induces
an antilinear isomorphism
$D_{w-j,\sigma}\to D_{j,\sigma}$, $b\mapsto\bar b$.
Following \cite[\S\S1.4.1--1.4.2]{KatoI}, for nonzero
$a\in D_{j,\sigma}$ and $b\in D_{w-j,\sigma}$, define
\[
\|a\otimes b^{-1}\|_{\mathrm{Hdg},\sigma}=|z|,
\qquad a=z\bar b.
\]
For $x\in L_j(P)$, set
\[
\|x\|_{j,v}=\|1\otimes x\|_{\mathrm{Hdg},\sigma}^{[F_v:\R]},
\]
with exponent $1$ at a real place and $2$ at a complex place
\cite[\S1.4.3]{KatoI}.
For $0\ne x\in L_j(P)$, we have $\|x\|_{j,v}=1$ at almost all places,
so we define
\begin{equation}\label{eq:basicheight}
H_j(P)=\prod_v\|x\|_{j,v}^{-1/2}.
\end{equation}
Since $\|ax\|_{j,v}=|a|_v\|x\|_{j,v}$ for $a\in F^\times$,
the product formula shows that $H_j(P)$ is independent of $x$.
We set $H_j(P)=1$ when the relevant graded pieces vanish.

For $V\in\Rep_\Q(G)$, fix a lattice $L_V\subset V$ with
$\mathcal K$-stable completion $\widehat L_V$, taking $L_0$ for $V_0$.
We use the transported lattice
\[
T(\Mf(V))=\widetilde\lambda(\widehat L_V).
\]
It is independent of the representative $\widetilde\lambda$ because
$\mathcal K$ preserves $\widehat L_V$, and it is Galois-stable
by the Galois invariance of $\lambda$.

\subsubsection{The canonical-type height}
The adjoint motive $A_{\mathcal M}$ is pure of weight zero, since $w_G$
acts trivially on $\mathfrak g$ by conjugation. The faithful representation
$G\hookrightarrow\operatorname{GL}(W)$ induces an injection
$\mathfrak g\hookrightarrow\End_\Q(W)$ on Lie algebras. Set
\[
 L_{\mathfrak g}:=\mathfrak g\cap\End_{\Z}(L).
\]
Its completion $\widehat L_{\mathfrak g}$ is stable under conjugation by $\mathcal K$.
We compute $H_j(A_{\mathcal M})$ using the transported lattice
\[
T(A_{\mathcal M})=\widetilde\lambda(\widehat L_{\mathfrak g}).
\]
\begin{definition}[Canonical-type height]
Following \cite[\S3.6.1]{KatoII}, define
\begin{equation}\label{eq:canonicalheight}
 H_{\spadesuit,\Phi}(\mathcal M)
   :=\prod_{j<0}H_j(A_{\mathcal M})^{-1},
 \qquad
 h_{\spadesuit,\Phi,K}(\mathcal M)
   :=\frac{\log H_{\spadesuit,\Phi}(\mathcal M)}{[F:K]}.
\end{equation}
\end{definition}

For the fixed embedding $F\subset\C$, choose $\theta$ as above and let
$h\in D_\Phi$ be the corresponding Hodge point. By Betti--de Rham
comparison and \eqref{eq:tangent}, the complexification of
$A_{\mathcal M,\mathrm{dR}}/\Fil^0A_{\mathcal M,\mathrm{dR}}$ identifies
with $T_{D_\Phi,h}$. A polarization of the Hodge realization of
$A_{\mathcal M}$ pairs its $j$th and $(-j)$th graded pieces, giving
\[
 \C\otimes_F\bigotimes_{j<0}L_j(A_{\mathcal M})^{-1}
 \simeq \bigl(\det T_{D_\Phi,h}^*\bigr)^{\otimes2}.
\]
This square of the canonical line accounts for the factor $1/2$ in
\eqref{eq:basicheight}. Kato's logarithmic Hodge interpretation
\cite[\S\S3.8.9--3.8.10 and 3.8.13]{KatoII} motivates the analogy with $K_X+D$.
The canonical-type height uses the adjoint motive of the fixed group $G$;
changing the ambient Lie algebra can change the height.

\subsubsection{A positive height for the same type}
To define a height analogous to $h_A$ in \eqref{eq:classical},
we require curvature positivity on the horizontal tangent spaces
of $D_\Phi$ in \eqref{eq:tangent}.

Let $J=\{j:\gr^j_{\Fil_h}V_{0,\C}\ne0\}$ for any $h\in D_\Phi$.
For $j\in J$, let $\mathscr H_j$ be the line bundle on $D_\Phi$ with
fiber $\det_\C\gr^j_{\Fil_h}V_{0,\C}$ at $h$. Let $\kappa_j$ be the
curvature form of its Hodge metric, with Kato's sign convention
\cite[\S\S2.1.2 and 3.5.3]{KatoII}.

\begin{definition}[Positive height datum]
For a real tuple $\Lambda=(c_j)_{j\in J}$, define
\[
 H_\Lambda(\mathcal M)=\prod_{j\in J}H_j(M_0)^{c_j},
 \qquad
 h_{\Lambda,K}(\mathcal M)
      =\frac{\log H_\Lambda(\mathcal M)}{[F:K]}.
\]
We call $\Lambda$ \emph{positive} if
$\kappa_\Lambda:=\sum_jc_j\kappa_j$ is positive definite on
$T_{D_\Phi,\mathrm{hor}}$. In this case, we call $h_{\Lambda,K}$ a
\emph{positive height}.
\end{definition}

For heights formed from $V_0$, this is the curvature condition in
Kato's definition of ampleness \cite[\S3.6.9]{KatoII}.
Here positivity refers to horizontal curvature, not to
nonnegativity of $h_{\Lambda,K}$.
If $T_{D_\Phi,\mathrm{hor}}=0$, every real tuple $\Lambda$
satisfies the condition.
At the Hodge point of $\mathcal M$, the fiber of
$\mathscr H_j\otimes\mathscr H_{w_0-j}^{-1}$ identifies with
$\C\otimes_F L_j(M_0)$. Since $q_0$ pairs the two graded pieces into a
constant Tate line, $\kappa_{w_0-j}=-\kappa_j$. Thus this line bundle
has curvature $2\kappa_j$.
The standard choice $c_j=j$ gives Kato's height
$H_\star(M_0)=\prod_jH_j(M_0)^j$ \cite[\S1.4.5]{KatoI}.
Since $W$ is faithful and its Tate summand has constant Hodge filtration,
every nonzero horizontal tangent vector at $h$ induces a nonzero map
$\gr^j_{\Fil_h}V_{0,\C}\to\gr^{j-1}_{\Fil_h}V_{0,\C}$ for some $j$.
Thus this datum is positive by Griffiths' curvature theorem
\cite[Theorem~2.2.3.1 and \S3.6.9]{KatoII}.

For comparison with Koshikawa's pure-motive height
\cite[\S4]{Koshikawa2015}, see \cite[\S\S1.4.10--1.4.12]{KatoI}.
The canonical line has positive horizontal curvature
\cite[Theorem~2.3.3]{KatoII}.
The canonical-type height $h_{\spadesuit,\Phi,K}$ attached to the adjoint
representation corresponds to $h_{K_X+D}$. The positive height
$h_{\Lambda,K}$ attached to the chosen datum $\Lambda$ corresponds to $h_A$.

\subsection{Degeneration and the counting term}
We define an untruncated counting term from the integral realizations,
following the local lattice constructions of \cite[\S2.3.4]{KatoI}
and \cite[\S3.5.8]{KatoIIarXiv}, specialized to the principal level
\eqref{eq:level} with $L\supset nL$ and trivial scalar subgroup.
Kato later did not use the arithmetic functions $H_{\heartsuit,v}$ and
$H_{\heartsuit,S_F}$ in \cite[Remark~3.6.13]{KatoII}.
We use this construction to measure monodromy at the fixed level;
only the nonnegativity of $N^\heartsuit_{K,S}$ is needed for our main theorem.

Fix a finite place $v$ of $F$ with residue field $k_v$ of
characteristic $p$. For a uniformizer $\varpi_v$, normalize
$\operatorname{ord}_v(\varpi_v)=1$.
For each prime $\ell$, let $T_\ell$ be the $\ell$-component of
\[
\widetilde\lambda(\widehat L)\subset\Mf(W)_{\mathrm{et}},
\]
and put $V_\ell=\Q_\ell\otimes_{\Z_\ell}T_\ell$.
Write
$
\rho_{v,\ell}:\operatorname{Gal}(\overline{F_v}/F_v)
\rightarrow\operatorname{GL}(V_\ell)
$
for the local Galois action.

For $\ell\ne p$, choose a generator of $\Z_\ell(1)$ and use it to
write the tame character of inertia associated with $\varpi_v$ as
$t_{v,\ell}:I_{F_v}\to\Z_\ell$. The nilpotent monodromy operator
$N'_{v,\ell}$ on $V_\ell$ is characterized by
\[
 \rho_{v,\ell}(\sigma)=\exp\bigl(t_{v,\ell}(\sigma)N'_{v,\ell}\bigr)
\]
for $\sigma$ in an open subgroup of $I_{F_v}$ \cite[\S1.7.2]{KatoI}.

At $p$, let $W(k_v)$ be the ring of $p$-typical Witt vectors, put
$F_{v,0}=\operatorname{Frac}W(k_v)$, and let $F_{v,0,\mathrm{ur}}$
be its maximal unramified extension in $\overline{F_v}$.
Since $V_p$ is de Rham, it becomes semistable over a finite
extension of $F_v$. Let
\[
D_v:=D_{\mathrm{pst}}(F_v,V_p)
\]
be its potentially semistable period module, of dimension
$\dim_{\Q_p}V_p$ over $F_{v,0,\mathrm{ur}}$.
Writing $\overline{k_v}$ for the residue field of $\overline{F_v}$,
put
\[
\widehat D_v:=W(\overline{k_v})[1/p]
\otimes_{F_{v,0,\mathrm{ur}}}D_v.
\]
Kato's construction \cite[\S1.3.7]{KatoI} gives a
$W(\overline{k_v})$-lattice
\[
\mathcal D_v(T_p):=D_{\mathrm{pst},W(\overline{k_v})}(T_p)
\subset\widehat D_v.
\]

Let $N^{\mathrm{abs}}_{v,p}$ be the nilpotent monodromy operator on
$D_v$ induced by $-d/ds_p$ on Fontaine's semistable period ring
$B_{\mathrm{st}}=B_{\mathrm{cris}}[s_p]$.
Here $B_{\mathrm{cris}}$ is the crystalline period ring and $s_p$
is the logarithmic variable associated with compatible $p$-power
roots of $p$ \cite[\S2.1]{Koshikawa2015}.
Following Kato's monodromy convention \cite[\S1.7.2]{KatoI},
we normalize the operator by setting
\[
N'_{v,p}:=\operatorname{ord}_v(p)\,N^{\mathrm{abs}}_{v,p}
\]
on $D_v$, and extend it to $\widehat D_v$.

Let $I_v(\mathcal M)$ be the set of all $b\in\Q$ such that
$\exp(bN'_{v,\ell})$ preserves $T_\ell$ and acts trivially on
$T_\ell/nT_\ell$ for every $\ell\ne p$, while
$\exp(bN'_{v,p})$ preserves $\mathcal D_v(T_p)$ and acts trivially on
$\mathcal D_v(T_p)/n\mathcal D_v(T_p)$.
These exponentials are finite sums because the operators are nilpotent.
Kato's construction shows that $I_v(\mathcal M)$ is a nonzero
subgroup of $\Q$ \cite[\S2.3.4]{KatoI}.

Define
\begin{equation}\label{eq:exponent}
e_v(\mathcal M)=
\begin{cases}
a^{-1},&I_v(\mathcal M)=a\Z\quad(a\in\Q_{>0}),\\
0,&I_v(\mathcal M)\text{ is not cyclic}.
\end{cases}
\end{equation}
The second case includes $I_v(\mathcal M)=\Q$.
The invariant $e_v(\mathcal M)$ is independent of the choice of
$\widetilde\lambda$ representing $\lambda$ and of the tame generator.

\begin{definition}[Counting term]
Set
\[
H_{\heartsuit,S_F}(\mathcal M)
      :=\prod_{v\notin S_F}(\#k_v)^{e_v(\mathcal M)},
\qquad
N^\heartsuit_{K,S}(\mathcal M)
      :=\frac1{[F:K]}
        \sum_{v\notin S_F}e_v(\mathcal M)\log\#k_v.
\]
\end{definition}

At a place of good reduction, the $\ell$-adic realizations are
unramified for $\ell\ne p$ and the $p$-adic realization is crystalline,
so all $N'_{v,\ell}$ vanish. Hence $I_v(\mathcal M)=\Q$ and
$e_v(\mathcal M)=0$.
Since a motive has good reduction at almost all places, the counting
term is a finite nonnegative sum.

Kato expected $e_v(\mathcal M)$ to be integral and its vanishing to
imply good reduction \cite[\S2.3.4]{KatoI}. The latter implication
can fail at a place dividing the level.
For example, take $n=3$,
$F=\Q(\zeta_3)$, and $E:y^2+y=x^3$, where $\zeta_3$ is a primitive
cube root of unity.
Its $3$-division polynomial is $3x(x^3+1)$, with
$y$-coordinates $0,-1$ for $x=0$ and $\zeta_3,\zeta_3^2$ for $x^3=-1$.
Thus $E[3]\subset E(F)$.
At $v\mid3$, the equation has discriminant $-27$ of valuation $6$.
Since the valuation of the discriminant is invariant modulo $12$
under changes of Weierstrass equation, $E$ has bad reduction at $v$.
Since $j(E)=0$, the curve has potentially good reduction.
Thus all nilpotent monodromy operators vanish, giving
$e_v(\mathcal M_E)=0$.

\subsection{Discriminants and base change}
Let $E/F/K$ be finite extensions in $\overline K$.
For a pure motive $P/F$ and a level object $\mathcal M/F$, write
$P_E$ and $\mathcal M_E$ for their base changes, with the transported
integral structures and level.
We have
\[
 H_j(P_E)=H_j(P)^{[E:F]},
 \qquad
 H_{\heartsuit,S_E}(\mathcal M_E)
     =H_{\heartsuit,S_F}(\mathcal M)^{[E:F]}.
\]
The first formula is \cite[Proposition~1.4.4 (4)]{KatoI}.
For the second, compare \cite[Proposition~3.5.14 (1)]{KatoIIarXiv}.

Let $w\mid v$ be finite places of $E/F$, and put
$e=e(w/v)$ and $p=\operatorname{char}(k_v)$.
Using a common algebraic closure of $F_v\subset E_w$,
restriction of the Galois action leaves $T_\ell$ unchanged.
Kato's construction identifies
$\mathcal D_w(T_p)$ with $\mathcal D_v(T_p)$
\cite[\S1.3.7]{KatoI}.
For $\ell\ne p$, the tame characters, using the same generator of
$\Z_\ell(1)$, satisfy $t_{v,\ell}|_{I_{E_w}}=e\,t_{w,\ell}$,
since a uniformizer at $v$ is a unit times the $e$th power of one at $w$.
For $\ell=p$, the operator $N^{\mathrm{abs}}_{v,p}$ is unchanged and
$\operatorname{ord}_w(p)=e\operatorname{ord}_v(p)$.
Thus $N'_{w,\ell}=eN'_{v,\ell}$ for every prime $\ell$.
It follows that $b\in I_w(\mathcal M_E)$ if and only if
$eb\in I_v(\mathcal M)$, so
\[
 I_w(\mathcal M_E)=e^{-1}I_v(\mathcal M),
 \qquad
 e_w(\mathcal M_E)=e\,e_v(\mathcal M),
\]
including the noncyclic case. Since
$\sum_{w\mid v}e(w/v)\log\#k_w=[E:F]\log\#k_v$,
summing over $v\notin S_F$ gives the second formula. In particular,
\begin{equation}\label{eq:basechange}
 \begin{split}
 h_{\spadesuit,\Phi,K}(\mathcal M_E)
      &=h_{\spadesuit,\Phi,K}(\mathcal M),\\
 h_{\Lambda,K}(\mathcal M_E)&=h_{\Lambda,K}(\mathcal M),\\
 N^{\heartsuit}_{K,S}(\mathcal M_E)
      &=N^{\heartsuit}_{K,S}(\mathcal M).
 \end{split}
\end{equation}
The discriminant is not base-change invariant.
Indeed, the discriminant tower formula
\cite[\S24, (117)]{Vojta} gives
\[
 d_K(E)=d_K(F)+
 \frac{\log N_{F/\Q}\mathfrak d_{E/F}}{[E:K]}\geq d_K(F),
\]
where $\mathfrak d_{E/F}$ denotes the relative discriminant ideal.

\section{The motivic Vojta inequality and our results}
We formulate a bounded-degree version of Kato's motivic Vojta
question for pure Andr\'e motives.
We verify the inequality for several classes of motives.

\subsection{The motivic Vojta inequality}
\label{subsec:motivic-vojta}
Let $\mathscr G_{\Phi,r}(K)$ be the disjoint union of the motivic generic classes in $\mathscr M_{\Phi,\mathcal K}(F)$ over the embedded fields $K\subseteq F\subset\overline K$ with $[F:K]\leq r$.

The main subject of this paper is the following bounded-degree,
untruncated variant of Kato's question for pure Andr\'e motives \cite[\S4.2.3, Question~4 (1)]{KatoII},
using the modifications discussed in \cite[Remark~4.2.2.1]{KatoIIarXiv}.

\begin{question}[A pure-motive variant of Kato's question]
\label{conj:main}
Fix $K,S,r,\Phi$, the lattice and level, and a positive
height datum $\Lambda$. For which such data does the following bound hold?
For every $\varepsilon>0$ there is a real
constant $C_\varepsilon$, depending only on these fixed data and
$\varepsilon$, such that
\begin{equation}\label{eq:main}
 h_{\spadesuit,\Phi,K}(\mathcal M)
 \leq d_K(F)+N^{\heartsuit}_{K,S}(\mathcal M)
                  +\varepsilon h_{\Lambda,K}(\mathcal M)+C_\varepsilon
\end{equation}
for every $(F,\mathcal M)\in\mathscr G_{\Phi,r}(K)$.
\end{question}

\begin{remark}[Known cases of the classical conjecture]
\label{rem:classical-known-cases}
For $r=1$, \eqref{eq:classical} is known for smooth projective
curves of genus $0$, $1$, and at least $2$ by Roth's theorem,
Lang's approximation theorem, and Faltings's finiteness theorem,
respectively.
It also holds for projective space with hyperplanes in general
position by Schmidt's Subspace Theorem, and for abelian varieties
with an ample reduced normal crossings divisor by Faltings's
approximation theorem
\cite[Theorem~8.10 and \S15, pp.~158--159]{Vojta}.
\end{remark}

The correspondence with the classical inequality \eqref{eq:classical} is
\[
\begin{array}{c|c}
\text{classical term}&\text{Motivic term}\\ \hline
x/K(x)&\mathcal M/F\\
h_{K_X+D}(x)&h_{\spadesuit,\Phi,K}(\mathcal M)\\
d_K(K(x))&d_K(F)\\
N_S(D,x)&N^{\heartsuit}_{K,S}(\mathcal M)\\
h_A(x)&h_{\Lambda,K}(\mathcal M).
\end{array}
\]
Without further geometric comparison results, the classical Vojta
conjecture does not directly imply the bound in Question~\ref{conj:main}.

We now compare our formulation with Kato's in several respects.

\smallskip
\noindent\emph{Kato's question.}
In \cite[\S4.2.3]{KatoII}, Kato considers arithmetic and Hodge--Nevanlinna analogues of Vojta's conjecture in Questions~4 (1) and~4 (2), respectively.
We compare our inequality with the former.
Kato's question is formulated for arbitrary fixed data.
We ask instead for which fixed data \eqref{eq:main} holds uniformly on
$\mathscr G_{\Phi,r}(K)$.
We write $K,S,S_F$ for Kato's $F_0,S_0,S$, respectively.

\smallskip
\noindent\emph{Group-relative data.}
We fix $G$, its tensor data, and the integral level, and define
the canonical-type height using $\Mf(\Lie G)$.
Our objects use the absolute level condition of \cite[\S1.3.3]{KatoII}.
If a level object over $K$ is chosen, Kato identifies the absolute
parameter sets over finite extensions of $K$ with the relative ones
for $\mathcal G=G$, after conjugating the level subgroup
\cite[\S1.3.7(2)]{KatoII}.
His relative framework also allows a normal subgroup
$\mathcal G\triangleleft G$ and a fixed quotient motive on
$\Rep(G/\mathcal G)$ \cite[\S\S1.3.4--1.3.5]{KatoII}.
We use the unconditional Tannakian category of pure Andr\'e
motives, whereas Kato uses Jannsen's mixed motives
\cite[\S1.3.1]{KatoII}.
For matching realizations and integral lattices, our basic heights
agree with Kato's, while our canonical-type height is the logarithm
of his group-relative height divided by $[F:K]$
\cite[\S\S3.2, 3.6.1, and 3.6.7 (2)]{KatoII}.

\smallskip
\noindent\emph{Discriminant normalization.}
We use the full relative discriminant $d_K(F)$, as in the classical Vojta conjecture.
Kato instead uses the discriminant outside $S$
\cite[\S4.2.1]{KatoII}:
\[
d_{K,S}(F):=\frac{\log D_S(F/K)}{[F:K]},
\]
where $D_S(F/K)$ is the norm of the part of the different
outside $S_F$.

Both terms have coefficient $1$, and for fixed $K,S,r$,
there is a constant $B_{K,S,r}$ such that
\begin{equation}\label{eq:discriminant-outside-S}
0\leq d_K(F)-d_{K,S}(F)\leq B_{K,S,r}
\qquad([F:K]\leq r).
\end{equation}
Indeed, only places above $S$ contribute to the difference,
with uniformly bounded local different exponents.
Thus the two terms are interchangeable up to a uniform
additive constant.

\smallskip
\noindent\emph{The untruncated variant.}
Kato's Question~4(1) uses the truncated counting term
\[
N^{(1)}_{K,S}(\mathcal M)
:=\frac1{[F:K]}
\sum_{v\in\Sigma_{S_F}(\mathcal M)}\log\#k_v,
\]
where $\Sigma_{S_F}(\mathcal M)$ is the set of bad-reduction
places outside $S_F$.
For matching pure data, its logarithmic form is
\begin{equation}\label{eq:kato-truncated}
h_{\spadesuit,\Phi,K}(\mathcal M)
\leq d_{K,S}(F)+N^{(1)}_{K,S}(\mathcal M)
+h_{\Lambda,K}(\mathcal M)+C,
\end{equation}
with a uniform constant $C$ and no degree bound.

Our term $N^\heartsuit_{K,S}$ weights each place by the monodromy
multiplicity $e_v(\mathcal M)$, paralleling the untruncated counting
term in \cite[Conjecture~25.1 (b)]{Vojta}.
Kato discusses a degree bound and replacing the truncated counting
term with $H_{\heartsuit,S_F}$ \cite[Remark~4.2.2.1]{KatoIIarXiv}.
Combining these leads to \eqref{eq:main} for matching data,
by \eqref{eq:discriminant-outside-S} and rescaling $\Lambda$.
The definition of $H_{\heartsuit,S_F}$ requires $S_F$ to contain
the bad-reduction places of the base motive \cite[\S3.5.8]{KatoIIarXiv}.
Our inequality allows any fixed $S$ containing the archimedean places.
If $e_v(\mathcal M)\geq1$ at every bad-reduction place outside $S_F$, then
$N^{(1)}_{K,S}\leq N^\heartsuit_{K,S}$, so the untruncated inequality is weaker than the truncated one.
The full level-$3$ example in \S2.3 shows that this lower bound
can fail at a bad-reduction place outside $S_F$.

\smallskip
\noindent\emph{Positivity and motivic genericity.}
Our positive height data specialize Kato's ample data to the pure
case, using only the factors $H_j(\Mf(V_0))$.
His mixed extension-height factors are trivial here, so we can take
the required $t(w,d)>0$ without changing the height
\cite[\S\S3.2.3 and 3.6.9]{KatoII}.

Our condition $G_{\mathcal M}^0=G^0$ is Kato's Tannakian
genericity condition, applied in the Andr\'e category
\cite[\S1.3.13]{KatoII}.
Hodge genericity implies this condition, but the converse is
not assumed.

\smallskip
\noindent\emph{The role of $\varepsilon$.}
For every $\varepsilon>0$, the datum $\varepsilon\Lambda$ remains
positive, and
$h_{\varepsilon\Lambda,K}=\varepsilon h_{\Lambda,K}$.
Thus quantifying over all positive height data with coefficient $1$
is equivalent to allowing every $\varepsilon>0$ for each fixed datum
\cite[Remark~4.2.3.1]{KatoII}.
We keep $\varepsilon$ explicit, as in the classical Vojta conjecture.

\begin{remark}[The exceptional-set limitation]
Motivic genericity does not replace Vojta's ex\-cep\-tion\-al-set
condition: a proper subvariety of a geometric parameter space
may still contain objects with $G_{\mathcal M}^0=G^0$.
\end{remark}

\subsection{Strong form with an arbitrary additive constant}
A stronger form of the bound in Question~\ref{conj:main} allows the additive
constant to be prescribed arbitrarily:
\begin{equation}\label{eq:strong}
 \text{for every }\varepsilon>0\text{ and }C\in\R,
 \quad\eqref{eq:main}\text{ with }C_\varepsilon=C
 \text{ holds outside a finite set.}
\end{equation}
Question~\ref{conj:main} asks only for one constant for each $\varepsilon$.
\begin{proposition}\label{prop:northcott}
Assume the following fixed-level finiteness property: for every
$B\in\R$, there are finitely many pairs $(F_i,\mathcal M_i)$ with
$K\subseteq F_i\subset\overline K$ and
$\mathcal M_i\in\mathscr M_{\Phi,\mathcal K}(F_i)$ such that every
$(F,\mathcal M)\in\mathscr G_{\Phi,r}(K)$ with
$h_{\Lambda,K}(\mathcal M)\leq B$ is isomorphic to a base change of
some $\mathcal M_i$ along $F_i\subseteq F$, including the level and
integral structures. Under this assumption,
the bound in Question~\ref{conj:main} is equivalent to \eqref{eq:strong}.
\end{proposition}

\begin{proof}
Assume the bound in Question~\ref{conj:main}. Fix $\varepsilon>0$ and $C\in\R$.
We show that only finitely many objects fail \eqref{eq:main} with
$C_\varepsilon=C$.
Applying the bound with $\varepsilon/2$ gives a constant $C_0$.
Then Any such failure satisfies
\[
 h_{\Lambda,K}(\mathcal M)<\frac{2(C_0-C)}{\varepsilon}.
\]
By the finiteness assumption, these objects are base changes of
finitely many $\mathcal M_i$.
For each such base change, \eqref{eq:basechange} and
$N^\heartsuit_{K,S}\geq0$ give
\[
 d_K(F)<h_{\spadesuit,\Phi,K}(\mathcal M_i)
                  -\varepsilon h_{\Lambda,K}(\mathcal M_i)-C.
\]
Thus $d_K(F)$ is uniformly bounded.
Since $[F:K]\leq r$, Hermite--Minkowski gives only finitely many
fields $F$, including their embeddings in $\overline K$.
Hence only finitely many such base changes occur, proving
\eqref{eq:strong}.

Conversely, the finitely many exceptions in \eqref{eq:strong}
can be absorbed into a larger uniform constant, giving
the bound in Question~\ref{conj:main}.
\end{proof}

\subsection{A single Hodge type}

In the following subsections, we verify \eqref{eq:main} for several
classes of motivically generic objects.
For each class, we fix the tensor type, lattices, principal level
$n\geq3$, and positive datum $\Lambda$.
The constants may depend on the fixed parameters of each class.

We begin with the simplest case, where both height terms vanish
and the inequality follows from nonnegativity.

\begin{proposition}
\label{prop:examples-diagonal}
Suppose that $V_0$ has only Hodge type $(a,a)$. Then
\[
 h_{\spadesuit,\Phi,K}(\mathcal M)
 =h_{\Lambda,K}(\mathcal M)=0.
\]
Consequently, \eqref{eq:main} holds with $C_\varepsilon=0$.
\end{proposition}

\begin{proof}
Let $h:\mathbb S\to G_\R$ be the Hodge homomorphism associated with
$\mathcal M$.
It acts by scalars on both $G$-stable summands
of the faithful representation $W=V_0\oplus\Q_\eta$.
Thus $h$ is central in $G_\R$, so $A_{\mathcal M}$ has only
Hodge type $(0,0)$; hence \eqref{eq:canonicalheight} gives
$h_{\spadesuit,\Phi,K}(\mathcal M)=0$.
Since $J=\{a\}$ and $L_a(M_0)$ is a line tensored with its metric
dual, $H_a(M_0)=1$ and $h_{\Lambda,K}(\mathcal M)=0$.
The assertion follows from $d_K(F),N^\heartsuit_{K,S}\geq0$.
\end{proof}

\begin{example}[Artin--Tate motives]
The proposition applies to $M_0=A(-a)$ with $A$ an Artin motive;
our convention gives this motive Hodge type $(a,a)$.
In fact, $A$ is split over $F$. Put $d=\dim_\Q A_B$.
After removing the Tate twist, the transported adelic lattice is the
completion of a lattice in $A_B$. The principal level condition makes
the Galois action trivial modulo $n$ on both $M_0$ and the Tate line.
Thus the finite
Galois image of $A$ lies in the principal congruence subgroup of
$\operatorname{GL}_d(\Z)$. This subgroup is torsion-free for $n\geq3$,
so $A\simeq\Q(0)^{\oplus d}$.
Then all local monodromy operators vanish, so
\eqref{eq:exponent} gives $e_v(\mathcal M)=0$ for every finite
place $v$. Hence $N^\heartsuit_{K,S}(\mathcal M)=0$.

For these split examples, take
$G=\Gm$, $\eta=\operatorname{id}$, $w_G(t)=t^2$, and
$V_0=(\Q_\eta^{\otimes a})^{\oplus d}$, with a polarizing pairing.
The Tate object remains part of the generated category, so their
connected motivic group is $\Gm$.
Thus they are generic for this type.
In particular, the pure components
\[
 h^{2i}(\mathbb P^m_F)\simeq\Q(-i),\qquad
 h^0(C)\simeq\Q(0),\qquad h^2(C)\simeq\Q(-1)
\]
provide examples, where $C/F$ is a smooth projective geometrically connected curve of genus zero.
The full motive of $\mathbb P^m$ for $m>0$ has several weights and is not an example of the fixed pure type.
\end{example}

\subsection{Height formulas for abelian varieties}
\label{subsec:abelian-heights}
We recall the height formulas used in both the CM and non-CM cases.
For an abelian variety $A/F$, we use the natural integral realization
\[
T(h^1(A))=\prod_\ell
H^1_{\mathrm{et}}(A_{\overline F},\Z_\ell)
\]
and the usual integral Tate line. Put
\[
b_K(A):=\frac{\log H_1(h^1(A))}{[F:K]},
\qquad
\beta_\Lambda:=c_1-c_0.
\]
We relate $b_K(A)$ to the absolute stable Faltings height
$h_{\mathrm{Fal}}^{\mathrm{st}}(A)$, with the normalization compatible
with Kato's height.

Choose a finite extension $F'/F$ over which $A$ has semistable
reduction. By \cite[Example~1.4.20]{KatoI},
\[
\log H_1(h^1(A_{F'}))
=[F':\Q]h_{\mathrm{Fal}}^{\mathrm{st}}(A).
\]
By \cite[Proposition~1.4.4 (4)]{KatoI}, this gives
\[
\log H_1(h^1(A))
=\frac{\log H_1(h^1(A_{F'}))}{[F':F]}
=[F:\Q]h_{\mathrm{Fal}}^{\mathrm{st}}(A).
\]
Dividing by $[F:K]$, we obtain
\begin{equation}\label{eq:b-faltings-height}
b_K(A)=[K:\Q]h_{\mathrm{Fal}}^{\mathrm{st}}(A).
\end{equation}

Next, \cite[Proposition~1.4.4 (3)]{KatoI} gives
\[
H_0(h^1(A))=H_1(h^1(A))^{-1}.
\]
Since $h^1(A)$ has only Hodge types $(1,0)$ and $(0,1)$,
an object $\mathcal M$ with $M_0=h^1(A)$ and the above integral
structure satisfies
\[
H_\Lambda(\mathcal M)
=H_0(h^1(A))^{c_0}H_1(h^1(A))^{c_1}
=H_1(h^1(A))^{c_1-c_0}.
\]
Taking normalized logarithms gives
\begin{equation}\label{eq:examples-abelian-height}
h_{\Lambda,K}(\mathcal M)
=(c_1-c_0)b_K(A)
=\beta_\Lambda b_K(A).
\end{equation}

\subsection{CM abelian varieties}
We consider CM abelian varieties of a fixed tensor type with $G^0$ a torus.
The canonical-type height and the degeneration term vanish, while the
positive height need not.

For principally polarized $A/F$, a principal level-$n$ structure
with these integral structures requires $A[n]\subset A(F)$,
and the Weil pairing then gives $\mu_n\subset F$.

We use the fixed integral tensor data to control this remaining height term.

\begin{proposition}[CM examples of a fixed integral tensor type]
\label{prop:examples-cm}
Fix a dimension $g$ and a weight-one type $\Phi$ with
$\dim V_0=2g$ and $G^0$ a torus. For the motivic generic natural level
objects of this fixed type arising from principally polarized CM
abelian varieties $A/F$, inequality \eqref{eq:main} holds over every
finite extension $F/K$, with a constant independent of $[F:K]$.
This holds for every positive datum $\Lambda$, including
$\beta_\Lambda<0$. In this type,
\[
h_{\spadesuit,\Phi,K}(\mathcal M)=0,
\qquad
N^\heartsuit_{K,S}(\mathcal M)=0.
\]
\end{proposition}

\begin{proof}
Put $T=G^0$, and let $h:\mathbb S\to G_\R$ be the Hodge homomorphism associated with $\mathcal M$.
Since $\mathbb S$ is connected, $h$ factors through $T_\R$.
Since $T$ is a torus, it acts trivially on
$\Lie(G)=\Lie(T)$, so $A_{\mathcal M}$ has only Hodge type $(0,0)$.
Thus \eqref{eq:canonicalheight} gives
$h_{\spadesuit,\Phi,K}(\mathcal M)=0$, and \eqref{eq:tangent} gives
$T_{D_\Phi,\mathrm{hor}}=0$. Hence every real tuple $\Lambda$ is positive.

A CM abelian variety has potentially good reduction at every finite place \cite[Theorem~6 (a)]{SerreTateGoodReduction}.
Thus the local base-change formulas give vanishing of all monodromy operators.
By \eqref{eq:exponent}, $e_v(\mathcal M)=0$ at every finite place,
so $N^\heartsuit_{K,S}(\mathcal M)=0$.

It remains to bound $b_K(A)$.
Conjugation by $T(\R)$ fixes $h$, so
$D_\Phi=G(\R)\cdot h$ is finite.
Borel's class-number theorem \cite[Theorem~5.1]{BorelFiniteness}
gives
\[
\#\bigl(G(\Q)\backslash G(\Af)/\mathcal K\bigr)<\infty.
\]
Thus
$
\mathcal Y_{\Phi,\mathcal K}
:=G(\Q)\backslash
\bigl(D_\Phi\times G(\Af)/\mathcal K\bigr)
$
is finite.
Every object determines a class in $\mathcal Y_{\Phi,\mathcal K}$
by Betti--\'etale comparison and $\theta$
\cite[\S\S1.2.2 and 1.2.5--1.2.6]{KatoII}.
This class determines the integral weight-one Hodge structure of $A$:
its Hodge filtration comes from $D_\Phi$, and its integral lattice
from $\widehat L_0$ transported by $\widetilde\lambda$ and viewed
on the Betti side. Such a Hodge structure determines the complex
abelian variety up to isomorphism, so only finitely many geometric
isomorphism classes of $A$ occur.
Stable Faltings height is invariant under geometric isomorphism, so
$
|b_K(A)|\leq B_{\Phi,L,\mathcal K,K}
$
for a constant depending only on the fixed data.
By \eqref{eq:examples-abelian-height},
$h_{\Lambda,K}(\mathcal M)=\beta_\Lambda b_K(A)$.
Since $d_K(F)\geq0$, taking
$
C_\varepsilon
=\varepsilon|\beta_\Lambda|B_{\Phi,L,\mathcal K,K}
$
proves \eqref{eq:main}.
\end{proof}

\begin{remark}[Scope]
\label{rem:examples-cm-degree}
For a fixed tensor type $\Phi$ with $G^0$ a torus and a fixed integral level,
the proposition holds uniformly over all finite extensions $F/K$.
It applies only to CM abelian varieties compatible with these data:
arbitrary commensurable lattices are excluded, and we do not claim
that one fixed tensor type with $G^0$ a torus covers all CM types of dimension $g$.
\end{remark}

\subsection{Height formulas for elliptic curves}
\label{subsec:elliptic-heights}
To verify \eqref{eq:main} for several non-CM elliptic curves, we first collect the height and level formulas for general elliptic curves.
Let $E/F$ be any elliptic curve, with $F/K$ finite.

We compute the canonical-type height.
Let $\Phi_{\mathrm{ell}}$ denote the usual elliptic type:
\[
G=\operatorname{GL}_2,\qquad
V_0=\Q^2,\qquad
\eta=\det,\qquad
w_G(t)=tI_2,
\]
with the standard alternating pairing and $L_0=\Z^2$.
On $W=\Q^2\oplus\Q_{\det}$, $\mathfrak g=\mathfrak{gl}_2$ acts by
$X\mapsto(X,\operatorname{tr}X)$.
Thus $L_{\mathfrak g}=\mathfrak g\cap\End_\Z(L)$ is identified with
$\End_\Z(L_0)$.
For an object $\mathcal M$ of this type with $M_0=h^1(E)$
and the natural integral structure, the adjoint representation gives
\[
A_{\mathcal M}=\Mf(\mathfrak g)\simeq\operatorname{End}(h^1(E)),
\]
with the natural endomorphism lattice.
The only negative Hodge graded piece of $A_{\mathcal M,\mathrm{dR}}$ is
\[
\gr^{-1}_{\Fil}A_{\mathcal M,\mathrm{dR}}
\simeq
\operatorname{Hom}_F
\bigl(\gr^1_{\Fil}h^1(E)_{\mathrm{dR}},
      \gr^0_{\Fil}h^1(E)_{\mathrm{dR}}\bigr).
\]
By \cite[\S1.4.14]{KatoI}, the natural isomorphism
\[
L_{-1}(\operatorname{End}(h^1(E)))
\simeq
L_0(h^1(E))\otimes L_1(h^1(E))^{-1}
\simeq
L_1(h^1(E))^{\otimes(-2)}
\]
preserves the metrics at every place. Hence
\[
H_{-1}(\operatorname{End}(h^1(E)))
=H_1(h^1(E))^{-2}.
\]
Since only the index $-1$ contributes to
\eqref{eq:canonicalheight}, we obtain
$h_{\spadesuit,\Phi_{\mathrm{ell}},K}(\mathcal M)=2b_K(E)$.

The pairing $q_0$ identifies $\gr^0_{\Fil}V_{0,\C}$ with the metric
dual of $\gr^1_{\Fil}V_{0,\C}$ tensored with the constant Tate line.
Hence $\kappa_0=-\kappa_1$, so
$\kappa_\Lambda=(c_1-c_0)\kappa_1=\beta_\Lambda\kappa_1$.
As $\kappa_1$ is positive definite on the elliptic horizontal tangent
space \cite[Theorem~2.2.3.1]{KatoII}, positivity of $\Lambda$ implies
$\beta_\Lambda>0$.
Thus we conclude that
\begin{equation}\label{eq:examples-elliptic-height}
h_{\spadesuit,\Phi_{\mathrm{ell}},K}(\mathcal M)=2b_K(E),
\qquad
\beta_\Lambda>0.
\end{equation}
This uses the full elliptic type $\Phi_{\mathrm{ell}}$, not the type with $G^0$ a torus used for CM curves.

We compare the stable Faltings height with the height of the
$j$-invariant. Write
\[
h_K(z)=[K:\Q]h(z),
\]
where $h$ is the absolute logarithmic Weil height on
$\mathbb P^1(\overline\Q)$. For a finite extension $L/K$, an element
$z\in L$, and a finite place $v$ of $L$, put
\[
m_v(z)=\max\{0,-\operatorname{ord}_v(z)\},\qquad q_v=\#k_v,
\]
with $m_v(0)=0$. Then
\begin{equation}\label{eq:examples-j-decomposition}
 h_K(z)
 =\frac1{[L:K]}
 \left(
   \sum_{v\text{ finite}}m_v(z)\log q_v
   +\sum_{\sigma:L\hookrightarrow\C}
       \log\max\{1,|\sigma(z)|\}
 \right).
\end{equation}

Our Faltings norm \cite[\S\S1.4.1--1.4.3]{KatoI} is
$\pi^{-1/2}$ times the norm in \cite[(10)]{Pazuki2019}.
Thus our stable Faltings height is larger by $\frac12\log\pi$.
Silverman's comparison \cite[(50)]{Pazuki2019} gives
\[
\frac1{12}h(j(E))
-\frac12\log(1+h(j(E)))-O(1)
\leq h_{\mathrm{Fal}}^{\mathrm{st}}(E)
\leq \frac1{12}h(j(E))+O(1).
\]
Since
$
\frac{x}{12}-\frac12\log(1+x)
$
is bounded below for $x\geq0$, and
$
\frac12\log(1+x)\leq\frac{x}{24}+O(1),
$
there are constants $C_K,C'_K\geq1$ such that, for every elliptic curve
over a finite extension of $K$,
\begin{equation}\label{eq:examples-height-comparison}
-C_K\leq b_K(E)\leq\frac1{12}h_K(j(E))+C_K,\qquad
h_K(j(E))\leq24b_K(E)+C'_K.
\end{equation}

The field $F(E[n])$ defines full level $n$, and the faithful action
of $\operatorname{Gal}(F(E[n])/F)$ on $E[n]$ gives
\begin{equation}\label{eq:examples-full-level-degree}
[F(E[n]):F]\leq|\operatorname{GL}_2(\Z/n\Z)|.
\end{equation}

We reduce \eqref{eq:main} to a numerical estimate for the stable
Faltings height.

For an elliptic curve $E/F$ with $E[n]\subset E(F)$, a
\emph{natural full level-$n$ object} is an object arising from
$h^1(E)$ with the geometric integral structures
and a principal level-$n$ orbit induced by a choice of basis of $E[n]$.
An inequality holds \emph{uniformly} on a class if a single constant
$C_\varepsilon$, depending only on the fixed data and additional
parameters of the class, works for all its curves, fields, and
associated level objects.

\begin{lemma}[A numerical criterion]
\label{lem:examples-numerical}
Fix a number field $K$, a finite set $S$ of places of $K$ containing
the archimedean places, an integer $n\geq3$, and a positive height
datum $\Lambda$ for $\Phi_{\mathrm{ell}}$.
Let $\mathscr E$ be a class of pairs $(F,E)$, where $F/K$ is finite and
$E/F$ is an elliptic curve with $E[n]\subset E(F)$.
Assume that, for every $\delta>0$, there is a real constant $C_\delta$
such that
\begin{equation}\label{eq:examples-numerical}
(2-\delta)b_K(E)\leq d_K(F)+C_\delta
\end{equation}
for all $(F,E)\in\mathscr E$.
Then, for every $\varepsilon>0$, \eqref{eq:main} holds uniformly for all
natural full level-$n$ objects of type $\Phi_{\mathrm{ell}}$ attached to
the curves in $\mathscr E$.

In particular, the assumption is satisfied if there is a real constant
$C$, independent of $(F,E)$, such that
\begin{equation}\label{eq:examples-strong-numerical}
2b_K(E)\leq d_K(F)+C
\end{equation}
for all $(F,E)\in\mathscr E$.

\end{lemma}

\begin{proof}
Fix $\varepsilon>0$ and put $\delta=\varepsilon\beta_\Lambda>0$.
By \eqref{eq:examples-abelian-height}, \eqref{eq:examples-elliptic-height}, and \eqref{eq:examples-numerical},
every such natural full level-$n$ object
$\mathcal M_E$ satisfies
\[
h_{\spadesuit,\Phi_{\mathrm{ell}},K}(\mathcal M_E)
-\varepsilon h_{\Lambda,K}(\mathcal M_E)
=(2-\delta)b_K(E)
\leq d_K(F)+C_\delta.
\]
Since $N^\heartsuit_{K,S}(\mathcal M_E)\geq0$, this proves \eqref{eq:main} uniformly.

Finally, if \eqref{eq:examples-strong-numerical} holds, then
\eqref{eq:examples-height-comparison} gives $b_K(E)\geq-C_K$, and hence
\[
(2-\delta)b_K(E)
\leq d_K(F)+C+\delta C_K
\qquad(\delta>0).
\]
Thus the assumption holds with $C_\delta=C+\delta C_K$.
\end{proof}

\begin{remark}[The truncated inequality]
Under \eqref{eq:examples-strong-numerical}, we have
\[
h_{\spadesuit,\Phi_{\mathrm{ell}},K}(\mathcal M_E)
\leq d_K(F)+C,
\]
without a counting term or a positive-height term.
If $[F:K]\leq r$, \eqref{eq:discriminant-outside-S} and
\eqref{eq:examples-numerical} imply Kato's truncated inequality
\eqref{eq:kato-truncated}, even with $\varepsilon h_{\Lambda,K}$
in place of $h_{\Lambda,K}$, since $N^{(1)}_{K,S}\geq0$.
\end{remark}

\begin{remark}[Transfer along bounded-degree isogenies]
\label{rem:examples-isogeny-transfer}
If $E'/F$ admits a geometric isogeny from $E$ of degree at most $M$,
then
\[
|b_K(E')-b_K(E)|
\leq\frac{[K:\Q]}2\log M
\]
by \cite[(13)]{Pazuki2019}, applied after an extension defining the
isogeny. Hence \eqref{eq:examples-numerical} for $E$ gives
\[
(2-\delta)b_K(E')
\leq d_K(F)+C_\delta
      +|2-\delta|\frac{[K:\Q]}2\log M.
\]
Thus \eqref{eq:examples-numerical} holds uniformly for $E'$ when
$M$ is fixed, and Lemma~\ref{lem:examples-numerical} applies whenever
$E'[n]\subset E'(F)$.
\end{remark}

\subsection{Non-CM elliptic curves with bounded stable height}
\label{subsec:bounded-heights}
We apply the numerical criterion to non-CM elliptic curves of bounded
stable Faltings height.
These curves are motivically generic for $\Phi_{\mathrm{ell}}$,
since their Mumford--Tate group is $\operatorname{GL}_2$.

\begin{proposition}[Bounded stable height]
\label{prop:examples-bounded-elliptic}
\label{prop:examples-finite-elliptic}
\label{prop:examples-fixed-isogeny}
Fix $B\in\R$. Inequality \eqref{eq:main} holds uniformly for the
natural full level-$n$ objects of type $\Phi_{\mathrm{ell}}$ arising
from non-CM elliptic curves $E/F$ over arbitrary finite extensions
$F/K$ with $b_K(E)\leq B$. The constant is independent of $F$ and
$[F:K]$.

In particular, the inequality holds uniformly for these natural
full level-$n$ objects in each of the following cases:
\begin{enumerate}
\item \emph{(Finite geometric classes)}
$j(E)$ belongs to a fixed finite set $\mathcal J\subset\overline{\Q}$,
with no bound on $[F:K]$.
\item \emph{(Fixed geometric isogeny classes)}
$[F:K]\leq r$ and $E$ is geometrically isogenous to one of a fixed
finite collection of non-CM elliptic curves over $K$.
\end{enumerate}
The constants may depend on $\mathcal J$ in \textup{(1)},
and on the fixed curves and $r$ in \textup{(2)}.
\end{proposition}

\begin{proof}
Since $b_K(E)\leq B$ and $d_K(F)\geq0$, we have
\[
2b_K(E)\leq2B\leq d_K(F)+2B.
\]
Thus \eqref{eq:examples-strong-numerical} holds uniformly, and
Lemma~\ref{lem:examples-numerical} proves the bounded-height assertion.

\textup{(1)} the $j$-invariant determines the geometric isomorphism class of an elliptic curve.
Hence only finitely many stable Faltings heights occur,
so \textup{(1)} follows.

\textup{(2)} choose a curve $E_0/K$ in the fixed collection
which is geometrically isogenous to $E$.
By \cite[Th\'eor\`eme~1.4]{GaudronRemond2014}, there is a geometric
isogeny $E_0\to E$ whose degree is bounded in terms of $[F:\Q]$ and
$h_{\mathrm{Fal}}^{\mathrm{st}}(E_0)$. Since $[F:K]\leq r$, its degree
is uniformly bounded.
Thus the height estimate in Remark~\ref{rem:examples-isogeny-transfer} bounds $b_K(E)$ uniformly in terms of $E_0$, $K$, and $r$.
Taking the maximum over the fixed finite collection of curves gives a uniform height bound,
so \textup{(2)} follows.
\end{proof}

If $[F:K]\leq r$, a uniform bound on $b_K(E)$ also bounds
$h_K(j(E))$ by \eqref{eq:examples-height-comparison}.
Since $[K(j(E)):K]\leq r$, Northcott's theorem gives only
finitely many geometric isomorphism classes.

\subsection{A \texorpdfstring{$j$}{j}-invariant criterion for non-CM elliptic curves}
\label{subsec:j-criterion}

We give a sufficient condition for \eqref{eq:examples-strong-numerical} using the finite poles of $j(E)$.
It applies to non-CM elliptic curves with unbounded stable Faltings height.
We begin with a local divisibility property that will yield a lower bound for the discriminant.

\begin{lemma}[Local divisibility under trivial projective torsion]
\label{lem:examples-local-divisibility}
Let $k/\Q_p$ be a finite extension, with normalized valuation
$\operatorname{ord}_k:k^\times\to\Z$, and let $n\geq2$.
Let $E/k$ be an elliptic curve with
$\operatorname{ord}_k(j(E))<0$.
If the Galois action on $E[n]$ has scalar image, then
\[
n\mid-\operatorname{ord}_k(j(E)).
\]
No assumption that $p\nmid n$ is needed.
\end{lemma}

\begin{proof}
Put $G_k=\operatorname{Gal}(\overline k/k)$.
Tate uniformization gives a parameter $q\in k^\times$ with
$0<|q|_k<1$ whose Tate curve $E_q/k$ satisfies
\[
j(E_q)=j(E),\qquad
\operatorname{ord}_k(q)=-\operatorname{ord}_k(j(E));
\]
see \cite[Chapter~V, \S\S3--5]{SilvermanAdvanced}.
Since $j(E)\neq0,1728$, the curve $E$ is a quadratic twist of $E_q$
\cite[Chapter~X, \S5]{SilvermanAEC}.
After identifying $E[n]$ with $E_q[n]$ via the twist, their Galois
actions differ by a character with values in $\{\pm1\}$.
Hence $G_k$ acts on $E_q[n]$ by scalars in $(\Z/n\Z)^\times$.

Use the $G_k$-equivariant identification
$E_q(\overline k)\simeq\overline k^\times/q^\Z$.
For $[x]\in E_q[n]$ with $x^n=q^m$, the map
$\pi([x])=m\bmod n$ gives a $G_k$-equivariant exact sequence
\[
0\longrightarrow\mu_n
\longrightarrow E_q[n]
\xrightarrow{\pi}\Z/n\Z
\longrightarrow0
\]
\cite[``Points of finite order'']{TateUniformization}.
Since $G_k$ acts trivially on $\Z/n\Z$,
every such scalar equals $1$, and hence $G_k$ acts trivially on $E_q[n]$.
Thus every point of $\pi^{-1}(1)$ is $G_k$-fixed.
The map $x\mapsto[x]$ gives a $G_k$-equivariant bijection from the
set of $n$th roots of $q$ to $\pi^{-1}(1)$.
Hence $q=\alpha^n$ for some $\alpha\in k^\times$, and
\[
-\operatorname{ord}_k(j(E))
=\operatorname{ord}_k(q)
=n\operatorname{ord}_k(\alpha).
\]
This proves the assertion.

\end{proof}

For $m\geq1$, put
\[
c_n(m)=1-\frac{\gcd(n,m)}{n}.
\]
For a finite extension $L/K$ and $z\in L$, define
\begin{equation}\label{eq:examples-R}
R_{n,L/K}(z)
=\frac1{[L:K]}
\sum_{\substack{v\text{ finite}\\m_v(z)>0}}
c_n(m_v(z))\log q_v.
\end{equation}
This sum includes all finite poles, including those whose residue
characteristic divides $n$.

We now pass from the local divisibility condition to a lower bound for the discriminant.

\begin{proposition}[A discriminant estimate under trivial projective \(n\)-torsion]
\label{prop:examples-projective-level}
Let $K\subseteq L\subseteq F$ be number fields, and let $n\geq2$.
Let $E/F$ be an elliptic curve with $j=j(E)\in L$.
Assume that the representation
\[
\rho_{E,n}:\operatorname{Gal}(\overline F/F)
\rightarrow\operatorname{GL}_2(\Z/n\Z)
\]
has scalar image, equivalently, its image in
$\operatorname{GL}_2(\Z/n\Z)/(\Z/n\Z)^\times$ is trivial.
Then
\begin{equation}\label{eq:examples-projective-discriminant}
d_K(F)\geq d_K(L)+R_{n,L/K}(j).
\end{equation}
The assumption holds if $E[n]\subset E(F)$.
Under the assumption, the estimate applies with
$L=K(j)$, without requiring $E$ to descend to $K(j)$.
\end{proposition}

\begin{proof}
Fix a finite pole $v$ of $j$, and put
$m=m_v(j)=-\operatorname{ord}_v(j)>0$.
For every $w\mid v$, Lemma~\ref{lem:examples-local-divisibility},
applied to $E_{F_w}$, gives
\[
n\mid-\operatorname{ord}_w(j)=e(w/v)m.
\]
Equivalently,
\begin{equation}\label{eq:examples-projective-divisibility}
a_v:=\frac{n}{\gcd(n,m)}\mid e(w/v).
\end{equation}

Write the different of $F_w/L_v$ as
$\mathfrak D_{F_w/L_v}=\mathfrak m_w^{\delta_{w/v}}$,
where $\mathfrak m_w$ is the maximal ideal of $\mathcal O_{F_w}$.
The bound $\delta_{w/v}\geq e(w/v)-1$ \cite[Chapter~III, \S6, Proposition~13]{SerreLocalFields}
and \eqref{eq:examples-projective-divisibility} imply
\[
\delta_{w/v}
\geq e(w/v)-1
\geq\left(1-\frac1{a_v}\right)e(w/v).
\]
Thus,
\[
\operatorname{ord}_v(\mathfrak d_{F/L})
=\sum_{w\mid v}f(w/v)\delta_{w/v}
\geq c_n(m)\sum_{w\mid v}e(w/v)f(w/v)
=c_n(m)[F:L].
\]

By the ideal norm formula,
\[
\log N_{L/\Q}\mathfrak d_{F/L}
=\sum_{v\text{ finite}}
\operatorname{ord}_v(\mathfrak d_{F/L})\log q_v.
\]
All summands are nonnegative, so the bounds at the finite poles of $j$ give
\[
\frac{\log N_{L/\Q}\mathfrak d_{F/L}}{[F:K]}
\geq
\frac{[F:L]}{[F:K]}
\sum_{\substack{v\text{ finite}\\m_v(j)>0}}
c_n(m_v(j))\log q_v
=R_{n,L/K}(j).
\]

The discriminant tower formula
\[
d_K(F)
=d_K(L)
+\frac{\log N_{L/\Q}\mathfrak d_{F/L}}{[F:K]}
\]
proves \eqref{eq:examples-projective-discriminant}.
\end{proof}

\begin{remark}[Relation to level covers]
A related ramification argument appears in
\cite[\S4.4, Proposition~4.3]{AbramovichVarillyAlvarado2017},
where level covers give lower bounds for boundary intersection
multiplicities.
\end{remark}

For an integer $d\geq2$, put
\begin{equation}\label{eq:examples-rational-R}
\mathcal R_n(d):=R_{n,\Q/\Q}(-1/d)
=\sum_{p\mid d}c_n(\operatorname{ord}_p(d))\log p.
\end{equation}
For every integer $a$ with $\gcd(a,d)=1$, the finite poles of
$a/d$ have the same orders as those of $-1/d$, so
\[
R_{n,\Q/\Q}(a/d)=\mathcal R_n(d).
\]

The following criterion relates the height of a rational
$j$-invariant to the contribution from its finite poles
and plays a key role in the proof of our main result.

\begin{proposition}[A criterion for rational $j$-invariants]
\label{prop:examples-rational-j}
Fix a number field $K$, an integer $n\geq3$, and $B\geq0$.
Let $a,d$ be integers with $d\geq2$ and $\gcd(a,d)=1$.
Let $E/F$ be an elliptic curve over a finite extension $F/K$ with
$j(E)=a/d$, and assume that $\rho_{E,n}$ has scalar image.
\begin{enumerate}
\item
If
\begin{equation}\label{eq:examples-rational-height-condition}
\log\max\{|a|,d\}\leq6\mathcal R_n(d)+B,
\end{equation}
then
$2b_K(E)\leq d_K(F)+O_{K,B}(1)$.
\item
Fix $J_0\geq1$. If $|a|\leq J_0d$ and
\begin{equation}\label{eq:examples-rational-denominator}
\log d\leq6\mathcal R_n(d)+B,
\end{equation}
then
$2b_K(E)\leq d_K(F)+O_{K,J_0,B}(1).$
\end{enumerate}
Both bounds are uniform in $a$, $d$, $E$, and $F$.
In either case, if $E[n]\subset E(F)$, its natural full level-$n$
objects satisfy \eqref{eq:main} uniformly for every $\varepsilon>0$.
\end{proposition}

\begin{proof}
For \textup{(1)}, since $j(E)=a/d$ is in lowest terms, we have
\[
h_K(j(E))=[K:\Q]\log\max\{|a|,d\}.
\]
Combining this with \eqref{eq:examples-height-comparison} and
\eqref{eq:examples-rational-height-condition}, we obtain
\begin{equation}\label{ineq:examples-reciprocal-height-2}
2b_K(E)
\leq [K:\Q]\mathcal R_n(d)
   +\frac{[K:\Q]}6B+2C_K.
\end{equation}

Since $\rho_{E,n}$ has scalar image,
Proposition~\ref{prop:examples-projective-level}, applied with base and
intermediate field $\Q$, gives
\[
\mathcal R_n(d)
=R_{n,\Q/\Q}(j(E))
\leq d_{\Q}(F)
=\frac{d_K(F)+\log|D_K|}{[K:\Q]}.
\]
Combining this with \eqref{ineq:examples-reciprocal-height-2} gives
\[
2b_K(E)
\leq d_K(F)+\log|D_K|
+\frac{[K:\Q]}6B+2C_K.
\]

For \textup{(2)}, since $|a|\leq J_0d$, we have
\[
\log\max\{|a|,d\}
\leq\log d+\log J_0
\leq6\mathcal R_n(d)+B+\log J_0.
\]
Part~\textup{(1)}, with $B+\log J_0$ in place of $B$, gives the
stated bound.

When $E[n]\subset E(F)$, the final assertion follows from
Lemma~\ref{lem:examples-numerical}.
\end{proof}
Part~\textup{(2)} includes the case of bounded numerator.

\begin{remark}[Cube-free denominators]
\label{rem:examples-cube-free}
Every cube-free integer $d\geq2$ satisfies
\[
\log d\leq4\mathcal R_n(d),
\]
and hence \eqref{eq:examples-rational-denominator} with $B=0$.
Indeed, for each $p\mid d$, put
$m_p=\operatorname{ord}_p(d)\in\{1,2\}$.
Since $c_n(m)\geq m/4$ for $m=1,2$ and every $n\geq3$, we have
\[
\mathcal R_n(d)
=\sum_{p\mid d}c_n(m_p)\log p
\geq\frac14\sum_{p\mid d}m_p\log p
=\frac14\log d.
\]
\end{remark}

\begin{remark}[One-factor sufficient conditions]
\label{rem:examples-one-factor}
Fix $A\geq1$, and let $D\geq2$ have a positive divisor $Q$ such that
\[
\gcd(Q,D/Q)=1,\qquad D\leq A Q^2.
\]
If $Q$ is squarefree, or if $Q$ is cube-free and $n\neq4$, then
\eqref{eq:examples-rational-denominator} holds for the denominator $D$
with $B=\log A$.
Indeed, every prime $p\mid Q$ occurs in $D$ with the same exponent
$m_p=\operatorname{ord}_p(Q)$. Under either assumption, $m_p\leq2$
and $c_n(m_p)\geq m_p/3$, since
\[
c_n(1)=1-\frac1n\geq\frac13,
\qquad
c_n(2)=1-\frac{\gcd(n,2)}{n}\geq\frac23
\quad(n\neq4).
\]
Restricting the defining sum for $\mathcal R_n(D)$ to primes dividing
$Q$ gives
\begin{equation}\label{eq:examples-one-factor-R}
\begin{aligned}
6\mathcal R_n(D)
&\geq6\sum_{p\mid Q}c_n(m_p)\log p
\geq2\sum_{p\mid Q}m_p\log p\\
&=2\log Q.
\end{aligned}
\end{equation}
Moreover, the size assumption gives
\begin{equation}\label{eq:examples-one-factor-size}
\log D\leq2\log Q+\log A.
\end{equation}
Combining these inequalities proves the claim.

For \(n=4\), the cube-free one-factor condition is not sufficient in general.
For example, $A=1$, $Q=25$, and $D=400=2^4\cdot5^2$ satisfy
$\gcd(Q,D/Q)=1$ and $D\leq Q^2$, but
\[
6\mathcal R_4(D)=3\log5<\log400.
\]
Here $c_4(4)=0$ and $c_4(2)=\frac12<\frac23$.
The squarefree condition on \(Q\) remains sufficient.
\end{remark}

\subsection{Non-CM elliptic curves with unbounded heights}
\label{subsec:unbounded-heights}

We prove the main result, establishing
\eqref{eq:main} for natural level objects arising from infinitely many
geometrically distinct non-CM elliptic curves with unbounded stable
Faltings height and uniformly bounded full-level degree.

\begin{theorem}[Rational $j$-invariants and unbounded stable height]
\label{ex:examples-unbounded}
Fix a number field $K$, an integer $n\geq3$, and $J_0\geq1$.
Let $j=a/d\in\Q$ vary subject to
\[
a\in\Z,\qquad d\in\Z_{\geq2},\qquad
\gcd(a,d)=1,\qquad |a|\leq J_0d.
\]
Then
\begin{enumerate}
\item \emph{(Non-CM property and height estimate)}
Every elliptic curve $E/F$ with $F/K$ finite and $j(E)=j$ is non-CM
and satisfies
\begin{equation}\label{eq:examples-reciprocal-height}
b_K(E)=\frac{[K:\Q]}{12}\log d+O_{K,J_0}(1),
\end{equation}
uniformly in $a$, $d$, $E$, and $F$.

\item \emph{(The motivic Vojta inequality at level $n$)}
Fix $B\geq0$ and suppose that
\eqref{eq:examples-rational-denominator} holds.
Let $E/F$ be any elliptic curve as in \textup{(1)} for which
$\rho_{E,n}$ has scalar image. Then
\[
2b_K(E)\leq d_K(F)+C_{K,J_0,B},
\]
where $C_{K,J_0,B}$ is independent of $a$, $d$, $E$, and $F$.
If $E[n]\subset E(F)$, its natural full level-$n$ objects satisfy
\eqref{eq:main} uniformly for every $\varepsilon>0$.
By Remark~\ref{rem:examples-cube-free}, the conclusion holds whenever
$d$ is cube-free, with $B=0$.

\item \emph{(Bounded full-level degree and infinitely many classes)}
For each such $j$, let $E_j/K$ be an elliptic curve with $j(E_j)=j$.
For each quadratic twist $E_j^{(u)}$, with $u\in K^\times/(K^\times)^2$,
put $F_{j,u}=K(E_j^{(u)}[n])$. Then
\[
[F_{j,u}:K]\leq|\operatorname{GL}_2(\Z/n\Z)|.
\]
For any infinite set of such $j$ satisfying
\eqref{eq:examples-rational-denominator} with a fixed $B\geq0$,
the curves $E_j^{(u)}$ represent infinitely many geometric isomorphism
classes with unbounded stable Faltings height.
Their natural full level-$n$ objects over $F_{j,u}$ satisfy
\eqref{eq:main} uniformly in $j$ and $u$.
In particular, one may take $j=-1/p$ with $p$ prime.
\end{enumerate}
\end{theorem}

\begin{proof}
\smallskip
\noindent\emph{Proof of (1).}
Since $d\geq2$ and $\gcd(a,d)=1$, the rational number $j=a/d$ is not an algebraic integer.
Hence every $E$ is non-CM
\cite[Chapter~II, Theorem~6.1]{SilvermanAdvanced}.

We prove \eqref{eq:examples-reciprocal-height}.
Let $E/F$ be as in \textup{(1)}, and choose a finite
extension $F'/F$ over which $E$ is semistable.
Write $\mathfrak\Delta_{E/F'}$ for its minimal discriminant ideal.
Since $j(E)=a/d$ is in lowest terms, we have
\[
h(j(E))=\log\max\{|a|,d\},
\qquad
|\sigma(j(E))|=\frac{|a|}{d}\leq J_0
\]
for every embedding $\sigma:F'\hookrightarrow\C$.
Silverman's height formula \cite[(34)]{Pazuki2019} gives
\[
\frac{\log N_{F'/\Q}\mathfrak\Delta_{E/F'}}{[F':\Q]}
=h(j(E))-\log\max\{1,|a|/d\}
=\log d.
\]

Write $\Delta(\tau)=q\prod_{k\geq1}(1-q^k)^{24}$, where
$q=e^{2\pi i\tau}$.
Let
\[
\mathcal F=\{\tau\in\C:\operatorname{Im}\tau>0,\ 
|\operatorname{Re}\tau|\leq\tfrac12,\ |\tau|\geq1\}.
\]
For each $\sigma:F'\hookrightarrow\C$, choose $\tau_\sigma\in\mathcal F$
such that $(E_{F'})^\sigma(\C)\simeq\C/(\Z+\tau_\sigma\Z)$
as complex Lie groups.
Since $j(\tau)=q^{-1}+744+O(q)$ at the cusp and
$|j(\tau_\sigma)|=|\sigma(j(E))|\leq J_0$, the $\tau_\sigma$ lie in a
fixed compact subset of $\mathcal F$, depending only on $J_0$.
As $\Delta$ has no zeros there, we obtain
\[
\log\!\left(
(2\pi)^{12}|\Delta(\tau_\sigma)|
(\operatorname{Im}\tau_\sigma)^6
\right)=O_{J_0}(1)
\]
uniformly in $\sigma$, $a$, $d$, $E$, and $F'$.
Since $E_{F'}$ is semistable, Silverman's formula
\cite[(10), (35)]{Pazuki2019}, converted to our
Kato-compatible normalization, gives
\[
\begin{aligned}
h_{\mathrm{Fal}}^{\mathrm{st}}(E)
&=\frac{1}{12[F':\Q]}
\left(
\log N_{F'/\Q}\mathfrak\Delta_{E/F'}
-
\sum_{\sigma:F'\hookrightarrow\C}
\log\!\left(
(2\pi)^{12}|\Delta(\tau_\sigma)|
(\operatorname{Im}\tau_\sigma)^6
\right)
\right)
+\frac12\log\pi\\
&=\frac1{12}\log d+O_{J_0}(1).
\end{aligned}
\]
Combining this with \eqref{eq:b-faltings-height} gives
\[
b_K(E)
=[K:\Q]h_{\mathrm{Fal}}^{\mathrm{st}}(E)
=\frac{[K:\Q]}{12}\log d+O_{K,J_0}(1).
\]

\smallskip
\noindent\emph{Proof of (2).}
Proposition~\ref{prop:examples-rational-j} \textup{(2)} gives the
discriminant estimate when $\rho_{E,n}$ has scalar image,
and \eqref{eq:main} for natural full level-$n$ objects.
The cube-free case follows from Remark~\ref{rem:examples-cube-free}.

\smallskip
\noindent\emph{Proof of (3).}
Since $j\neq0,1728$, one may take
\[
E_j:\quad y^2=x^3-\frac{3j}{j-1728}x-\frac{2j}{j-1728},
\]
whose $j$-invariant is $j$.
For any choice of $E_j/K$ with $j(E_j)=j$ and any
$u\in K^\times/(K^\times)^2$, set
\[
F_{j,u}=K(E_j^{(u)}[n]).
\]
By definition, $E_j^{(u)}[n]\subset E_j^{(u)}(F_{j,u})$.
Applying \eqref{eq:examples-full-level-degree} to $E_j^{(u)}/K$ gives
\[
[F_{j,u}:K]\leq|\operatorname{GL}_2(\Z/n\Z)|,
\]
uniformly in $j$, $u$, and the choice of $E_j$.

Since only finitely many $j=a/d$ satisfy the conditions with bounded
$d$, we have $d\to\infty$ as $j$ runs through distinct values.
Hence Step~1 gives
\[
b_K(E_j^{(u)})=b_K(E_j)
=\frac{[K:\Q]}{12}\log d+O_{K,J_0}(1)\longrightarrow+\infty,
\]
uniformly in $u$.
Distinct $j$ give distinct geometric isomorphism classes.
If \eqref{eq:examples-rational-denominator} holds with a fixed $B$,
Step~2 gives \eqref{eq:main} over $F_{j,u}$ uniformly in $j$ and $u$.
Taking $j=-1/p$ for primes $p$ gives such an infinite collection
with $B=0$, since $p$ is cube-free and $|j|\leq1$.
Together with the degree bound, this proves the assertion;
in particular, $r\geq|\operatorname{GL}_2(\Z/n\Z)|$ suffices,
and for $n=3$ one may take $r\geq48$.

\end{proof}

\begin{remark}[A stronger estimate for cube-free denominators]
\label{rem:examples-cube-free-height}
Under the assumptions of Theorem~\ref{ex:examples-unbounded} \textup{(1)},
assume that $d$ is cube-free and $\rho_{E,n}$ has scalar image.
Then
\[
3b_K(E)\leq d_K(F)+O_{K,J_0}(1).
\]
Indeed, Remark~\ref{rem:examples-cube-free} and
\eqref{eq:examples-reciprocal-height} give
\[
3b_K(E)
=\frac{[K:\Q]}4\log d+O_{K,J_0}(1)
\leq[K:\Q]\mathcal R_n(d)+O_{K,J_0}(1).
\]
Proposition~\ref{prop:examples-projective-level}, applied with base and
intermediate field $\Q$, gives $d_{\Q}(F)\geq\mathcal R_n(d)$.
Since $d_K(F)=[K:\Q]d_{\Q}(F)-\log|D_K|$, the assertion follows.
\end{remark}

\begin{example}[An explicit Weierstrass family]
\label{ex:examples-weierstrass}
For each prime $p$, put $d_p=1728p+1$ and consider
\[
E_p:\ y^2=x^3-3d_px-2d_p^2.
\]
Here $c_4=144d_p$ and $\Delta=-1728d_p^3(d_p-1)$, so
$j=c_4^3/\Delta$ gives
\[
j(E_p)=-\frac{1728}{d_p-1}=-\frac1p.
\]
Thus the assumptions of Theorem~\ref{ex:examples-unbounded}
are satisfied with $a=-1$, $d=p$, $J_0=1$, and $B=0$,
since $p$ is cube-free.
Hence all conclusions of the theorem apply to this family.
\end{example}

We apply Theorem~\ref{ex:examples-unbounded} to a semistable
family. By Remark~\ref{rem:examples-one-factor}, we verify
\eqref{eq:examples-rational-denominator} using one factor of the denominator.

\begin{corollary}[A semistable family with unbounded stable height]
\label{cor:examples-semistable}
Fix a number field $K$ and an integer $n\geq3$.
For $t\in\Z_{\geq1}$, let
\[
\mathcal E_t:\ y^2+xy=x^3+t,\qquad
\mathcal D_t=t(432t+1).
\]
These curves are non-CM and semistable over $\Q$, with
\[
j(\mathcal E_t)=-1/\mathcal D_t.
\]
The stable-height estimate and full-level degree bound in
Theorem~\ref{ex:examples-unbounded} \textup{(1) and (3)}
apply to $\mathcal E_t$ and its quadratic twists, with
\[
a=-1,\qquad d=\mathcal D_t,\qquad J_0=1.
\]

Moreover, every natural full level-$n$ object attached to an $F$-form
of $\mathcal E_t$, with $F/K$ finite, satisfies \eqref{eq:main}
uniformly if at least one of the following conditions holds:
\begin{enumerate}
\item[\textup{(a)}] $\mathcal D_t$ is cube-free.
\item[\textup{(b)}] At least one of $t$ and $432t+1$ is squarefree.
\item[\textup{(c)}] $n\neq4$ and at least one of $t$ and $432t+1$ is cube-free.
\end{enumerate}

Thus prime values of $t$ yield infinitely many semistable non-CM
geometric isomorphism classes with unbounded stable Faltings height and
uniformly bounded full-level degree.
\end{corollary}

\begin{proof}
The standard Weierstrass formulas for $\mathcal E_t$ give
\[
c_4(\mathcal E_t)=1,\qquad
\Delta(\mathcal E_t)=-\mathcal D_t,\qquad
j(\mathcal E_t)=-1/\mathcal D_t.
\]
Since $c_4(\mathcal E_t)=1$, the equation is minimal at every prime.
At each prime $p\mid\mathcal D_t$, the reduction is multiplicative
of type $I_{\operatorname{ord}_p(\mathcal D_t)}$.
Thus $\mathcal E_t$ is semistable.

Since $\mathcal D_t\geq433$ and $|j(\mathcal E_t)|\leq1$,
Theorem~\ref{ex:examples-unbounded} \textup{(1) and (3)},
with $a=-1$, $d=\mathcal D_t$, and $J_0=1$,
gives the non-CM property, stable-height estimate, and
full-level degree bound.

It remains to prove \eqref{eq:examples-rational-denominator}
with $d=\mathcal D_t$ and $B=\log433$.
For \textup{(a)}, this follows from Remark~\ref{rem:examples-cube-free}.
For \textup{(b)}, take a squarefree $Q_t\in\{t,432t+1\}$;
for \textup{(c)}, take a cube-free $Q_t$ in the same set.
Since
\[
\gcd(t,432t+1)=1,\qquad
\mathcal D_t=t(432t+1)\leq433Q_t^2,
\]
Remark~\ref{rem:examples-one-factor}, with
$D=\mathcal D_t$, $Q=Q_t$, and $A=433$,
gives \eqref{eq:examples-rational-denominator} in both cases.
Theorem~\ref{ex:examples-unbounded} \textup{(2)} then gives
\eqref{eq:main} uniformly.

Finally, prime values of $t$ satisfy condition~\textup{(b)}.
Since $\mathcal D_t=t(432t+1)$ is strictly increasing and unbounded,
the $j$-invariants are pairwise distinct, and the stable-height
estimate gives unbounded heights. Together with the full-level
degree bound above, this proves the assertion.
\end{proof}

\begin{remark}[Earlier work on the family]
\label{rem:semistable-background}
Daniels studied the same family $\mathcal E_t$ and proved that
it is a Serre curve for every odd prime $t\ne7$
\cite[Theorem~1.6 and Example~8.2]{Daniels2015}.
\end{remark}

\subsection{Symmetric powers of the non-CM elliptic curves}
\label{subsec:symmetric-powers}

We extend \eqref{eq:main} from the preceding non-CM elliptic examples
to fixed symmetric powers of $h^1(E)$, with level induced from
full level $n$ on $E$.
The height comparisons reduce the proof to \eqref{eq:examples-numerical}.
Here symmetric powers are taken with the Tannakian
commutativity constraint of Andr\'e motives
\cite[\S4.3]{Andre}. Hence
\[
(\operatorname{Sym}^m h^1(E))_B
=\operatorname{Sym}^m H^1(E,\Q)
\]
has dimension $m+1$ for every $m\geq1$.

\begin{proposition}[Symmetric powers of the elliptic examples]
\label{prop:examples-symmetric}
Fix $m\geq1$, and let $E/F$ belong to one of the non-CM elliptic
subclasses for which \eqref{eq:examples-numerical} was established above.
Assume $E[n]\subset E(F)$.
Put
\[
P=\operatorname{Sym}^m h^1(E)
\]
and define
\[
G_m=\operatorname{im}\!\left(
\operatorname{GL}_2\rightarrow
\operatorname{GL}(\operatorname{Sym}^m\Q^2)\times\Gm,
\quad
g\mapsto(\operatorname{Sym}^m g,\det g)\right).
\]
Let $\Phi_m$ be the induced tensor type, with its induced polarization
and lattice
\[
L_{0,m}=\operatorname{Sym}^m\Z^2,
\qquad
L_m=L_{0,m}\oplus\Z e_\eta,
\]
and let $\mathcal K_m$ be the principal level subgroup
\eqref{eq:level} for $L_m$.

Then every level object on $P$ induced by a full level-$n$
structure on $E/F$ is motivically generic for $\Phi_m$ and
satisfies \eqref{eq:main} for every positive datum $\Lambda$ of this type.
The constants may also depend on $m$.
\end{proposition}

\begin{proof}
Let $\mathcal M_E$ be the given elliptic level object and
$\mathcal M_P$ its induced object.

\smallskip
\noindent\emph{Step 1: The induced group and level.}
Let $\varphi_m:\operatorname{GL}_2\to G_m$ be the homomorphism.
The kernel of $\operatorname{Sym}^m:\operatorname{GL}_2\to \operatorname{GL}(\operatorname{Sym}^m\Q^2)$
is the scalar subgroup $\mu_m$. For a scalar matrix $\zeta I_2$,
the second component of $\varphi_m$ is $\zeta^2$. Hence
\[
\ker(\varphi_m)=\mu_m\cap\mu_2
=\mu_{\gcd(m,2)}.
\]

Since this kernel is finite, the differential
\[
d\varphi_m:\mathfrak{gl}_2\xrightarrow{\sim}\Lie(G_m)
\]
is an isomorphism.
It is compatible with conjugation, so it identifies the adjoint
representations after pullback along $\varphi_m$.

The weight cocharacter of $G_m$ is
\[
w_m(t)=\varphi_m(tI_2)
=(t^m I_{\operatorname{Sym}^m\Q^2},t^2),
\]
and the Tate character is the second projection
$\eta_m:G_m\to\Gm$. Since $w_m(t)$ acts as $t^m$ on
$\operatorname{Sym}^m\Q^2$, $P$ has weight $m$.

The non-CM assumption gives the motivic group $\operatorname{GL}_2$ for $h^1(E)$ together with the Tate object.
The group generated by $P$ and $\Q(-1)$ is its image under $\varphi_m$, namely $G_m$.
Hence $P$ is motivically generic for $\Phi_m$.

Finally, an element acting trivially modulo $n$ on $\widehat{\Z}^2$ acts trivially modulo $n$ on $\operatorname{Sym}^m\widehat{\Z}^2$, and its determinant acts trivially modulo $n$ on the Tate line.
Thus the image of the elliptic principal level subgroup under $\varphi_m$ is contained in $\mathcal K_m$.
Therefore, a full level-$n$ structure on $E/F$ induces the required level orbit on $P$.

\smallskip
\noindent\emph{Step 2: The canonical-type height.}
Via $d\varphi_m$, we identify the adjoint motives and regard
$L_{\mathfrak g_m}$ and $L_{\mathfrak{gl}_2}$ as lattices in
$\mathfrak{gl}_2$, where $\mathfrak g_m=\Lie(G_m)$.
Each is full because it is the intersection of a rational Lie
algebra with the endomorphism lattice of its faithful representation.
Thus they are commensurable: there is an
integer $d_m\geq1$, depending only on $m$, such that
\begin{equation}\label{eq:examples-symmetric-lattice-comparison}
d_mL_{\mathfrak{gl}_2}
\subseteq L_{\mathfrak g_m}
\subseteq d_m^{-1}L_{\mathfrak{gl}_2}.
\end{equation}

Kato's integral construction is functorial for these fixed lattice
comparisons and for the corresponding filtered graded pieces
\cite[\S\S1.3.3--1.3.5]{KatoI}.
In particular, at a prime not dividing $d_m$, the two lattices agree.
The inclusions in \eqref{eq:examples-symmetric-lattice-comparison}
bound the local logarithmic metric differences
by $O_m(-\log|d_m|_v)$.
Hence their total finite-place difference is
\[
O_m\!\left(
\sum_{v\text{ finite}}-\log|d_m|_v
\right)
=
O_m([F:\Q]).
\]
Here we use $\sum_{v\text{ finite}}-\log|d_m|_v
=[F:\Q]\log d_m$.
At infinity, $d\varphi_m$ is an isomorphism of the adjoint real Hodge
structures, so it preserves Kato's metrics on their lines
$L_j$ \cite[\S1.4.1]{KatoI}. Thus there is no additional
archimedean error.
Thus, \eqref{eq:canonicalheight} and \eqref{eq:examples-elliptic-height} give
\[
h_{\spadesuit,\Phi_m,K}(\mathcal M_P)
=h_{\spadesuit,\Phi_{\mathrm{ell}},K}(\mathcal M_E)
+O_{m,K}(1)
=2b_K(E)+O_{m,K}(1).
\]

\smallskip
\noindent\emph{Step 3: The positive height.}
We first express the height lines of $P$ in terms of the elliptic
height line.
Put $A=\gr^1_{\Fil}h^1(E)_{\mathrm{dR}}$, $B=\gr^0_{\Fil}h^1(E)_{\mathrm{dR}}$.
Both are lines. For $0\leq j\leq m$,
\[
\gr^j_{\Fil}P_{\mathrm{dR}}
\simeq
A^{\otimes j}\otimes B^{\otimes(m-j)}.
\]
Applying this also with $m-j$ in place of $j$, we obtain
\[
\begin{aligned}
L_j(P)
&\simeq
\left(A^{\otimes j}\otimes B^{\otimes(m-j)}\right)
\otimes
\left(A^{\otimes(m-j)}\otimes B^{\otimes j}\right)^{-1}\\
&\simeq
(A\otimes B^{-1})^{\otimes(2j-m)}
=
L_1(h^1(E))^{\otimes(2j-m)}.
\end{aligned}
\]

Kato's integral construction is compatible with symmetric powers
\cite[\S1.3.4]{KatoI}. At a finite place $v\mid p$, the saturated
Hodge line in the rank-two $\mathcal O_{\C_v}$-lattice is a direct
summand. Choose a basis $u_1,u_0$ with $u_1$ spanning this line.
The class of $u_1^j u_0^{m-j}$ generates the $j$th graded lattice
of the symmetric power \cite[\S1.3.5]{KatoI}.
Thus these identifications preserve the finite-place metrics.
At infinity, both metrics and the symmetric-power
identifications are $\operatorname{GL}_2(\R)$-equivariant, so their
ratio is constant on each component of the elliptic Hodge domain.
These constants depend only on $m$. Summing over embeddings and
normalizing as in Step~2 gives, by the definition of $b_K(E)$,
\[
\frac{\log H_j(P)}{[F:K]}
=
(2j-m)b_K(E)+O_{m,K}(1).
\]

To compute $h_{\Lambda,K}$, we sum with the coefficients of
$\Lambda=(c_j)_{j=0}^m$ and obtain
\[
h_{\Lambda,K}(\mathcal M_P)
=
\alpha_{\Lambda,m}b_K(E)
+O_{m,K,\Lambda}(1),
\]
where
\[
\alpha_{\Lambda,m}
=
\sum_{j=0}^m c_j(2j-m).
\]

To apply the elliptic estimate in Step~4, we need
$\alpha_{\Lambda,m}>0$. We check this using the curvature
condition defining positivity of $\Lambda$.
Under the identification of horizontal tangent spaces induced by
$d\varphi_m$, the same graded-line description and
$\kappa_0(h^1(E))=-\kappa_1(h^1(E))$ give
\[
\kappa_j(P)=(2j-m)\kappa_1(h^1(E)),
\]
where $\kappa_j(P)$ is the curvature of $\det\gr^j_{\Fil}P$ and
$\kappa_1(h^1(E))$ is the curvature of the elliptic Hodge line.
Therefore
\[
\kappa_\Lambda
=
\alpha_{\Lambda,m}\kappa_1(h^1(E)).
\]
Since $\kappa_1(h^1(E))$ is positive on the elliptic horizontal direction,
the positivity of $\Lambda$ is equivalent to $\alpha_{\Lambda,m}>0$.
For the standard choice $c_j=j$,
\[
\alpha_{\Lambda,m}
=\sum_{j=0}^m j(2j-m)
=2\sum_{j=0}^m j^2
-m\sum_{j=0}^m j
=\frac{m(m+1)(m+2)}6.
\]

\smallskip
\noindent\emph{Step 4: Application of the elliptic estimate.}
The underlying elliptic subclass satisfies \eqref{eq:examples-numerical} uniformly.
Fix $\varepsilon>0$ and put
$\delta=\varepsilon\alpha_{\Lambda,m}>0$.
By Steps~2 and~3 and this numerical estimate,
\[
h_{\spadesuit,\Phi_m,K}(\mathcal M_P)
-\varepsilon h_{\Lambda,K}(\mathcal M_P)
=(2-\delta)b_K(E)+O_{m,K,\Lambda,\varepsilon}(1)
\leq d_K(F)+C_\delta+O_{m,K,\Lambda,\varepsilon}(1).
\]
Absorbing the bounded error into $C_\varepsilon$ and using
$N^\heartsuit_{K,S}(\mathcal M_P)\geq0$ gives \eqref{eq:main},
uniformly on the specified subclass.
\end{proof}

\begin{remark}[Scope]
The non-CM proofs establish \eqref{eq:main} without the nonnegative
term $N^\heartsuit_{K,S}$, so its local multiplicities
$e_v(\mathcal M)$ need not be computed.
For the unbounded-height elliptic families, ramification at finite
poles of $j(E)$ gives the discriminant bound, but the height estimate
is not automatic: $R_{n,L/K}$ ignores poles of order divisible by $n$
and the archimedean part of $h_K(j(E))$.
The symmetric-power result requires levels induced from the underlying elliptic curve.
\end{remark}

\subsection*{Acknowledgments}
The author would like to thank Professor Yu Yasufuku for helpful comments on an earlier draft,
especially regarding the discriminant term and the bounded-degree formulation
of Vojta's conjecture.
The author was supported by JSPS KAKENHI Grant number JP26K24498.
ChatGPT (OpenAI) was used in the preparation of this paper
to assist with mathematical exploration, checking arguments and computations,
and language editing.

\end{document}